\documentclass[a4paper,11pt]{amsart}

\usepackage{amssymb,amsmath,amsthm,mathrsfs,enumerate,graphicx, color}
\usepackage[pdfpagelabels,colorlinks,linkcolor=blue,citecolor=black,urlcolor=blue]{hyperref}
\usepackage{esint}
\usepackage{tikz}
\usetikzlibrary{arrows}
\newtheorem{thm}{Theorem}[section]

\newtheorem{cor}[thm]{Corollary}
\newtheorem{lem}[thm]{Lemma}
\newtheorem{prop}[thm]{Proposition}
\newtheorem{defn}[thm]{Definition}

\theoremstyle{definition}
\newtheorem{rem}[thm]{Remark}

\newcommand{\JJ}{\mathcal{J}}

\newcommand{\II}{\mathcal{I}}

\newcommand{\lesi}{\lesssim}
\newcommand{\x}{\mathbb{\times}}

\newcommand{\dx}{d\mu(x)}
\newcommand{\dy}{d\mu(y)}
\newcommand{\rx}{\rho(x)}

\newcommand{\ry}{\rho(y)}
\newcommand{\supp}{\operatorname{supp}}
\newcommand{\f}{\frac}

\newcommand{\Om}{\Omega}

\newcommand{\vc}{\infty}
\newcommand{\Rn}{\mathbb{R}^n}

\newcommand{\rad}{\rm rad}
\newcommand{\dm}{d\mu}

\title[Maximal function characterizations for  Hardy spaces and applications]{Maximal function characterizations for  Hardy spaces on spaces of homogeneous type with finite measure and applications}         
\author{The Anh Bui}
\address{Department of Mathematics, Macquarie University, NSW 2109,
Australia}
\email{the.bui@mq.edu.au, bt\_anh80@yahoo.com}
 
\author{Xuan Thinh Duong}
\address{Department of Mathematics, Macquarie University, NSW 2109,
	Australia}
\email{xuan.duong@mq.edu.au}

\author{Fu Ken Ly}
\address{The School of Mathematics and Statistics, The Faculty of Science, The University of Sydney, 
NSW 2006, Australia}
\email{ken.ly@sydney.edu.au}

\keywords{Hardy space, heat kernel, maximal function characterization, second order elliptic operator, Fourier--Bessel operator}
\subjclass[2010]{42B30, 42B35, 35K08, 35J25}

\begin{document}

\begin{abstract}
We prove nontangential and radial maximal function characterizations for Hardy spaces associated to a non-negative self-adjoint operator satisfying Gaussian estimates on a space of homogeneous type with finite measure. This not only addresses an open point in the literature, but also gives a complete answer to the question posed by Coifman and Weiss in the case of finite measure. We then apply our results to give maximal function characterizations for Hardy spaces associated to  second order elliptic operators with Neumann  and Dirichlet boundary conditions, Schr\"odinger operators with Dirichlet boundary conditions, and Fourier--Bessel operators.
\end{abstract}
\date{}

\maketitle

\tableofcontents

\section{Introduction}\label{sec: intro}

Let $(X,d, \mu)$ be a metric space endowed with a nonnegative Borel measure $\mu$ satisfying the following `doubling' condition: there exists a constant $C_1>0$ such that
\begin{equation}\label{doublingcondition}
\mu(B(x,2r))\leq C_1\mu(B(x,r))
\end{equation}
for all $x\in X$, $r>0$ and all balls $B(x,r):=\{y\in X: d(x,y)<r\}$. For the moment  $\mu(X)$ may be finite or infinite. 

It is not difficult to see that the condition \eqref{doublingcondition} implies that there exists a ``dimensional" constant $n\geq 0$ so that
\begin{equation}\label{doub2}
\mu(B(x,\lambda r))\leq C_2\lambda^n \mu(B(x,r))
\end{equation}
for all $x\in X, r>0$ and $\lambda\geq 1$, and
\begin{equation}\label{doub2s}
\mu(B(x, r))\leq C_3\mu(B(y,r))\Big(1+\f{d(x,y)}{r}\Big)^n
\end{equation}
for all $x,y\in X, r>0$.

 Assume also the existence of an operator $L$ that satisfies the following two conditions:
\begin{enumerate}
	\item[(A1)] $L$ is  a nonnegative self-adjoint operator on $L^2(X)$;
	\item[(A2)] $L$ generates a semigroup $\{e^{-tL}\}_{t>0}$ whose kernel $p_t(x,y)$ admits a Gaussian upper bound. That is, there exist two positive constants $C$ and  $c$ so that for all $x,y\in X$ and $t>0$,
	\begin{equation}
	\tag{GE}\label{GE}
	\displaystyle |p_t(x,y)|\leq \f{C}{\mu(B(x,\sqrt{t}))}\exp\Big(-\f{d(x,y)^2}{ct}\Big).
	\end{equation}
\end{enumerate}
Then for $0<p\le 1$ one can define three notions of Hardy spaces related to $L$. The first notion is through linear combinations of atoms that appropriately encode the cancellation inherent in $L$. The second and third notions  are $H^p_{L,\max}$ and $H^p_{L,\rad}$, which are defined via the non-tangential maximal function and the radial maximal function respectively. For the reader's  convenience we recall these notions below.

\begin{defn}[Atoms for $L$]\label{def: L-atom}
	Let $p\in (0,1]$ and $M\in \mathbb{N}$. A function $a$ supported in a ball $B$ is called a  $(p,M)_L$-atom if there exists a
	function $b\in {\mathcal D}(L^M)$ such that
	\begin{enumerate}[{\rm (i)}]
		\item  $a=L^M b$;
		\item $\supp L ^{k}b\subset B, \ k=0, 1, \dots, M$;
		\item $\|L^{k}b\|_{L^\vc(X)}\leq
		r_B^{2(M-k)}\mu(B)^{-\f{1}{p}},\ k=0,1,\dots,M$.
	\end{enumerate}
	In the particular case where $\mu(X)<\vc$, the constant function $[\mu(X)]^{-1/p}$ is also considered as an atom.
\end{defn}

\noindent Then the atomic Hardy space associated to the operator $L$ are defined as follows:
\begin{defn}[Atomic Hardy spaces for $L$]

	Given $p\in (0,1]$ and $M\in \mathbb{N}$, we  say that $f=\sum
	\lambda_ja_j$ is an atomic $(p,M)_L$-representation if
	$\{\lambda_j\}_{j=0}^\infty\in l^p$, each $a_j$ is a $(p,M)_L$-atom,
	and the sum converges in $L^2(X)$. The space $H^{p}_{L,at,M}(X)$ is then defined as the completion of
	\[
	\left\{f\in L^2(X):f \ \text{has an atomic
		$(p,M)_L$-representation}\right\},
	\]
	with the norm given by
	$$
	\|f\|_{H^{p}_{L,at,M}(X)}=\inf\left\{\left(\sum|\lambda_j|^p\right)^{1/p}:
	f=\sum \lambda_ja_j \ \text{is an atomic $(p,M)_L$-representation}\right\}.
	$$
\end{defn}

\noindent The maximal Hardy spaces associated to $L$ are defined as follows.

\begin{defn}[Maximal Hardy spaces for $L$]\label{defn-maximal Hardy spaces}

For $f\in L^2(X)$, we define the \textbf{non-tangential}  maximal function associated to $L$ of $f$ by
\[
f^*_{L}(x)=\sup_{t>0}\sup_{d(x,y)<t}|e^{-tL}f(y)|
\]
and the \textbf{radial} maximal function  by
\[
f^+_{L}(x)=\sup_{t>0}|e^{-tL}f(x)|.
\]
	Given $p\in (0,1]$,  the Hardy space $H^{p}_{L, {\rm max}}(X)$ is defined as the completion of
	$$
	\left\{f
	\in L^2(X): f^*_{L} \in L^p(X)\right\},
	$$
	with the norm given by
	$$
	\|f\|_{H^{p}_{L, {\rm max}}(X)}=\|f^*_{L}\|_{L^p(X)}.
	$$
	
	Similarly, the Hardy space $H^{p}_{L, {\rm rad}}(X)$ is defined as the completion of
	$$
	\left\{f
	\in L^2(X): f^+_{L} \in L^p(X)\right\},
	$$
	with the norm given by
	$$
	\|f\|_{H^{p}_{L,{\rm rad}}(X)}=\|f^+_{L}\|_{L^p(X)}.
	$$
\end{defn}
The theory of Hardy spaces associated to differential operators was initiated in \cite{ADM} and since then the theory has been studied intensively by many mathematicians. See for example \cite{DY, AMR, HM, HLMMY, JY} and the references therein. In this framework it is understood that the classical Hardy spaces $H^p(\Rn)$ can be viewed as the Hardy spaces associated to the Laplacian $-\Delta$.

A substantive problem in the theory of Hardy spaces is  to determine conditions for which the atomic and maximal notions coincide, and it this problem which is the focus of our paper.
More precisely we wish to answer the following question:

\textit{Question: Does the following equivalence hold:
\begin{equation}
\label{eq-maximal chracterization}
H^{p}_{L,at, M}(X) \equiv H^{p}_{L, {\rm max}}(X)\equiv H^{p}_{L, {\rm rad}}(X)
\end{equation}
for sufficiently large $M$?}

Before presenting our main result we highlight some history and known results related to  \eqref{eq-maximal chracterization}.
\begin{enumerate}[(i)]
\item In the Euclidean setting, when $L=-\Delta$, the Hardy spaces associated to $L$ and the classical Hardy spaces are identical. The classical Hardy spaces has its roots in complex function theory, and it was in that setting that the connection with the non-tangential maximal function was first elucidated \cite{BGS}. The role of maximal functions then took centre stage and was instrumental in the development of the real-variable theory beginning with the seminal work of Fefferman and Stein \cite{FS}. From that point onwards the theory developed rapidly and, through the efforts of \cite{Co,Latter,LU,Ca,CT},  the atomic characterization was added to the fold.

\smallskip

\item The notion of atoms enabled the extension of Hardy spaces from $\Rn$ to other structures \cite{CW}, and it was there that Coifman and Weiss introduced the concept of a space of homogeneous type. The viewpoint, as espoused in \cite{CW}, was to develop the theory on $X$ by starting with the notion of atomic Hardy spaces, which we shall denote by $H^p_{CW}(X)$ for $\f{n}{n+1}<p\le 1$ (see Definition \ref{def: Hp-CW} below). Under certain additional geometric assumptions, Coifman and
Weiss proved the radial maximal function characterization for $H^1_{CW}(X)$. They then proposed the following question: 

\textit{Question (Coifman-Weiss): Can one characterize the Hardy spaces $H^p_{CW}(X)$ by maximal functions for $p$ below 1?}

This question has been partly answered in the setting of Ahlfors $n$-regular metric measure spaces.
Recall that such spaces are spaces of homogeneous type
with $\mu(B(x,r))\sim r^n$ for all $x\in X$ and $r\in (0,2\, {\rm diam}(X))$. 
When $X$ is an Ahlfors $1$--regular metric measure space Uchiyama \cite{U} proved that the spaces $H^p_{CW}(X)$ can be characterized by radial maximal functions for $p<1$, but unfortunately the range of $p$ in \cite{U} is not optimal. The same result was obtained by \cite{MS} for the range $1/2<p\le 1$. A complete answer was given by \cite{YZ} but extra structural assumptions are needed -- namely a so called reverse-doubling condition on $\mu$ and that $\mu(X)=\infty$. To the best of our knowledge, the remaining case $\mu(X)<\infty$ is non-trivial and is still open.

\smallskip

\item In our setting, the theory of Hardy spaces arises from the fundamental observation that the classical Hardy spaces on $\Rn$ is intrinsically tied to the Laplacian $-\Delta$ and this observation allows the theory to be  generalized in another direction. The articles \cite{HM,HLMMY} give an account of this body of work and there one can also find partial answers to \eqref{eq-maximal chracterization}. The full equivalence was proved \cite{DKP,SY,SY2} but further assumptions were required in addition to (A1) and (A2). Reverse-doubling on $X$ and a regularity and markov condition on $L$ (see (A3) and (A4) below) was required in \cite{DKP}, while $\mu(X)=\infty$ was implicitly required in \cite{SY,SY2}. It is worth mentioning that the proofs in \cite{SY,SY2}, which are an adaptation of \cite{Ca}, does not work well in the case $\mu(X)<\infty$ and thus, in this situation, the problem is still open.
\end{enumerate}

This brings us to the first goal of the present article, which is to address the finite case in (iii) above. More precisely we shall prove
\begin{thm}\label{mainthm1}
	Let $\mu(X)<\infty$ and assume $L$ satisfies (A1) and (A2). Let $p\in (0,1]$, and $M>\f{n}{2}\big(\f{1}{p}-1\big)$. Then the Hardy spaces $H^{p}_{L,at,M}(X)$, $H^{p}_{L, {\rm max}}(X)$ and $H^{p}_{L, {\rm rad}}(X)$ coincide with equivalent norms.
\end{thm}
\noindent Let us explain the relevance of the condition $\mu(X)=\infty$ in \cite{SY,SY2}. The proofs there are rooted in decomposition of the product space $X\times(0,\infty)$, which we sketch here for the sake of convenience. For each $i\in \mathbb{Z}$ one defines the level set $O_i:=\{x\in X: \mathbf{M}f(x)>2^i\}$ where $\mathbf{M}$ is a certain maximal function that is lower-continuous, and the tent of ${O}_i$ through $\widehat{O}_i:=(x,t)\in X\times (0,\vc): B(x,4t)\subset O$. Then the space $X\times (0,\vc)$ can be decomposed as follows:
\begin{equation}\label{eq-Xvc}
X\times (0,\vc)= \bigcup_{i} \widehat{O}_i\backslash \widehat{O}_{i+1}.
\end{equation}
Unfortunately \eqref{eq-Xvc} fails in the case $X$ is bounded and this is the reason why the argument used in the case $\mu(X)=\vc$ is not applicable to the case $\mu(X)<\infty$. To overcome this obstacle, some new ideas are employed such as a new decomposition of  $X\times (0,\vc)$. It is  worth pointing out firstly  that our approach is also applicable for the case $\mu(X)=\vc$ and secondly, that  although our decomposition of the underlying product space $X\times (0,\vc)$ bears a resemblance to that in \cite{DKP}, the absence of both reverse-doubling on $X$ and the conditions (A3) and (A4) on $L$ requires some significant  innovations and improvements. The details can be found in Section 3.

By combining Theorem \ref{mainthm1}  with Theorem 1.2 of \cite{SY} we can now state the following, completing the picture in point (iii) above. 
\begin{cor}\label{cor1}
	Let $\mu(X)$ be finite or infinite and assume $L$ satisfies (A1) and (A2). Let $p\in (0,1]$, and $M>\f{n}{2}\big(\f{1}{p}-1\big)$. Then the Hardy spaces $H^{p}_{L,at,M}(X)$, $H^{p}_{L, {\rm max}}(X)$ and $H^{p}_{L, {\rm rad}}(X)$ coincide with equivalent norms. Due to this coincidence, we shall write $H^p_{L}(X)$ for any $H^{p}_{L,at,M}(X)$, $H^{p}_{L, {\rm max}}(X)$ and $H^{p}_{L, {\rm rad}}(X)$  for any such $p$ and $M$.
\end{cor}

\bigskip

Our second goal is to give the answer for the question in (ii) proposed by Coifman and Weiss under the presence of an operator $L$ when $\mu(X)$ is finite. We first recall the definition of the Hardy spaces $H^p_{CW}(X)$ on $X$. 

\begin{defn}[$p$-atoms]
Let $p\in(\f{n}{n+1},1]$. A function $a$ is called a $p$-atom associated to the ball $B$ if 
\begin{enumerate}[\upshape (i)]
	\item $\supp a\subset B$
	\item $\Vert a\Vert_{L^\infty(X)} \le \mu(B)^{-1/p}$
	\item $\displaystyle \int a(x)\,\dm(x)=0$
\end{enumerate}
When $\mu(X)<\infty$ then the constant function $\mu(X)^{-1/p}$ is also an atom. 
\end{defn}
\noindent To define the Hardy space $H^p_{CW}$ for $p$ below 1, we need to introduce the Lipschitz spaces $\mathfrak{L}_\alpha$. We say that the function $f$ is a member of $\mathfrak{L}_\alpha$ if there exists a constant $c>0$, such that
$$
|f(x)-f(y)|\leq c|B|^{\alpha}
$$
for all ball $B$ and $x, y\in B$, and the best constant $c$ can be taken to be the norm of $f$ and is denoted by $\|f\|_{\mathfrak{L}_\alpha}$.

\begin{defn}[Hardy spaces of Coifman and Weiss]\label{def: Hp-CW}
Let $\f{n}{n+1}<p\le 1$. We say that a function $f\in H^p_{CW}(X)$ if $f\in L^1(X)$ for $p=1$, or $f\in \mathfrak{L}_{1/p-1}^*$ for $p<1$, and there exists  a sequence $(\lambda_j)_{j\in \mathbb{N}}\in \ell^p$ and a sequence of $p$-atoms $(a_j)_{j\in \mathbb{N}}$ such that $f=\sum_{j}\lambda_ja_j$ in $L^1(X)$ for $p=1$, and $f=\sum_{j}\lambda_ja_j$ in $\mathfrak{L}_{1/p-1}^*$ for $p<1$. We set
\begin{align*}
\|f\|_{H^p_{CW}}=\inf\Big\{\Big(\sum_{j}|\lambda_j|^p\Big)^{1/p}: f=\sum_{j}\lambda_ja_j\Big\}.
\end{align*}
\end{defn}

We now consider the following two additional conditions for the operator $L$:
\begin{enumerate}
	\item[(A3)] There exists $\delta\in (0,1]$ such that for every $d(y,y')<\sqrt{t}/2$ and $0<t<{\rm diam}\,X$,
		\begin{align*}
	\displaystyle |p_t(x,y)-p_t(x,y')|\leq \Big(\f{d(y,y')}{\sqrt{t}}\Big)^\delta \f{C}{\mu(B(x,\sqrt{t}))}\exp\Big(-\f{d(x,y)^2}{ct}\Big) .
	\end{align*}
	\item[(A4)] For every $x\in X$ and $t>0$, we have
	\begin{align*}
	\displaystyle \int_X p_t(x,y)\dy =1.
	\end{align*}
\end{enumerate}
Then we have the following.
\begin{thm}\label{mainthm2}
Let $\mu(X)<\infty$ and assume the operator $L$ satisfies (A1), (A2), (A3) and (A4). Let $p\in (\f{n}{n+\delta},1]$ and $M>\f{n}{2}\big(\f{1}{p}-1\big)$. Then 
$H^p_{CW}(X)$, $H^p_{L,\max}(X)$ and $H^p_{L,\rad}(X)$ coincide with equivalent norms. 
\end{thm}
\noindent The equivalence in Theorem \ref{mainthm2} anwers the question proposed by Coifman and Weiss \cite{CW} mentioned in point (ii) above when $\mu(X)<\vc$.
Furthermore if (A3) is satisfied with $\delta=1$, then one obtains the optimal range $\f{n}{n+1}<p\le 1$.
\bigskip

The final aim of our article is to apply Theorems \ref{mainthm1} and \ref{mainthm2}  to certain differential operators on bounded/unbounded domains. We are able to prove the following new results:
\begin{enumerate}[{\rm (i)}]
	\item When $L$ is the second elliptic operator with Neumann boundary condition, we show that the Hardy spaces of extension or the Hard spaces of Coifman and Weiss  $H^p_{CW}(\Om)$ coincides with $H^p_{L,\max}(\Om)$ and $H^p_{L,\rad}(\Om)$. This not only extends the results  in \cite{AR} to $p<1$, but also furnishes a new result even for $p=1$with bounded $\Om$. See Theorem \ref{mainth-Neumann}.
	
	\item When $L$ is a second elliptic operator with Dirichlet boundary condition, we show that the Hardy spaces of Miyachi  $H^p_{Mi}(\Om)$ coincides with $H^p_{L,\max}(\Om)$ and $H^p_{L,\rad}(\Om)$. Our paper is the first to give  maximal function characterizations for the Hardy spaces of Miyachi  $H^p_{Mi}(\Om)$. Furthermore, in the particular case when $\Om$ is a strongly Lipschitz domain such that either $\Om$ is bounded or $
	\Om^c$ is unbounded, then we have $H^p_r(\Om)\equiv H^p_{L,\max}(\Om)\equiv H^p_{L,\rad}(\Om)$, improving the results of \cite{AR} for $p<1$. Here $H^p_r(\Om)$ is the Hardy spaces of restriction.  See Theorem \ref{mainth-Dirichlet}.
	
	\item When $L$ is a Schr\"odinger operators with Dirichlet boundary condition, we introduce a new version of Hardy space of Miyachi type  $H^p_{\rho}(\Om)$ and show that $H^p_r(\Om)\equiv H^p_{L,\max}(\Om)\equiv H^p_{L,\rad}(\Om)$. See Theorem \ref{mainth-SchrodingerDirichlet}.
	
	\item When $L$ is a Fourier--Bessel operator on $((0,1), dx)$ or $((0,1), x^{2\nu+1}dx)$, we show that the maximal Hardy spaces $H^p_{L,\max}(\Om)$ and $ H^p_{L,\rad}(\Om)$ enjoy certain atomic characterizations which extends the results in \cite{DPRS} to $p<1$ and a larger range of $\nu$. See Theorem \ref{Th1-Bessel 1} and Theorem \ref{Th1}.
 \end{enumerate}

We close this introduction with some remarks on the organization of the article. Section \ref{sect-prelim} collects  some useful estimates for the operator $L$ arising from (A1) and (A2). Section \ref{sect-mainproofs} contain the proofs of Theorems \ref{mainthm1} and \ref{mainthm2}, while the applications can be found in Section \ref{sect-app}.

\bigskip

\textbf{Notation.} As usual we use $C$ and $c$ to denote positive constants that are independent of the main parameters involved but may differ from line to line. The notation $A\lesi B$  denotes $A\leq CB$, and $A\sim B$ means that both $A\lesi B$ and $B\lesi A$ hold. We use $\fint_E f\dm=\f{1}{\mu(E)}\int_E f\dm$ to denote the average of $f$ over $E$.
We write $B(x,r)$ to denote the ball centred at $x$ with radius $r$. By a `ball $B$' we mean the ball $B(x_B, r_B)$ with some fixed centre $x_B$ and radius $r_B$. The annuli around a  given ball $B$ will be denoted by $S_j(B)=2^{j+1}B\backslash 2^jB$ for $j\ge 1$ and $S_0(B)=2B$ for $j=0$. 


\section{Some kernel  and maximal function estimates}\label{sect-prelim}

Let  $L$  satisfy (A1) and (A2). Denote by $E_L(\lambda)$ a spectral decomposition of $L$. Then by spectral theory, for any bounded Borel funtion $F:[0,\vc)\to \mathbb{C}$ we can define
$$
F(L)=\int_0^\vc F(\lambda)dE_L(\lambda)
$$
as a bounded operator on $L^2(X)$. It is well-known that the kernel $K_{\cos(t\sqrt{L})}$ of $\cos(t\sqrt{L})$ satisfies 
\begin{equation}\label{finitepropagation}
{\rm supp}\,K_{\cos(t\sqrt{L})}\subset \{(x,y)\in X\times X:
d(x,y)\leq t\}.
\end{equation}
See for example \cite{CS}.
We have the following useful lemmas.
\begin{lem}[\cite{HLMMY}]\label{lem:finite propagation}
	Let $\varphi\in C^\vc_0(\mathbb{R})$ be an even function with {\rm supp}\,$\varphi\subset (-1, 1)$ and $\int \varphi =2\pi$. Denote by $\Phi$ the Fourier transform of $\varphi$. Then for every $k\in \mathbb{N}$, the kernel $K_{(t\sqrt{L})^k\Phi(t\sqrt{L})}(\cdot, \cdot)$ of $(t\sqrt{L})^k\Phi(t\sqrt{L})$ satisfies 
	\begin{equation}\label{eq1-lemPsiL}
	\displaystyle
	{\rm supp}\,K_{(t\sqrt{L})^k\Phi(t\sqrt{L})}\subset \{(x,y)\in X\times X:
	d(x,y)\leq t\},
	\end{equation}
	and
	\begin{equation}\label{eq2-lemPsiL}
	|K_{(t\sqrt{L})^k\Phi(t\sqrt{L})}(x,y)|\leq \f{C}{\mu(B(x,t))}.
	\end{equation}
\end{lem}
\begin{lem}
	\label{lem1}
	\begin{enumerate}[{\rm (a)}]
		\item Let $\varphi\in \mathscr{S}(\mathbb{R})$ be an even function. Then for any $N>0$ there exists $C$ such that 
		\begin{equation}
		\label{eq1-lema1}
		|K_{\varphi(t\sqrt{L})}(x,y)|\leq \f{C}{\mu(B(x,t))+\mu(B(y,t))}\Big(1+\f{d(x,y)}{t}\Big)^{-n-N},
		\end{equation}
		for all $t>0$ and $x,y\in X$.
		\item Let $\varphi_1, \varphi_2\in \mathscr{S}(\mathbb{R})$ be even functions. Then for any $N>0$ there exists $C$ such that
		\begin{equation}
		\label{eq2-lema1}
		|K_{\varphi_1(t\sqrt{L})\varphi_2(s\sqrt{L})}(x,y)|\leq C\f{1}{\mu(B(x,t))+\mu(B(y,t))}\Big(1+\f{d(x,y)}{t}\Big)^{-n-N},
		\end{equation}
		for all $t\leq s<2t$ and $x,y\in X$.
		\item Let $\varphi_1, \varphi_2\in \mathscr{S}(\mathbb{R})$ be even functions with $\varphi^{(\nu)}_2(0)=0$ for $\nu=0,1,\ldots,2\ell$ for some $\ell\in\mathbb{Z}^+$. Then for any $N>0$ there exists $C$ such that
		\begin{equation}
		\label{eq3-lema1}
		|K_{\varphi_1(t\sqrt{L})\varphi_2(s\sqrt{L})}(x,y)|\leq C\Big(\f{s}{t}\Big)^{2\ell} \f{1}{\mu(B(x,t))+\mu(B(y,t))}\Big(1+\f{d(x,y)}{t}\Big)^{-n-N},
		\end{equation}
		for all $t\geq s>0$ and $x,y\in X$.
	\end{enumerate}
\end{lem}
\begin{proof}
	\noindent (a) The estimate \eqref{eq1-lema1} was proved in \cite[Lemma 2.3]{CD} in the particular case $X=\mathbb{R}^n$ but the proof is still valid in the spaces of homogeneous type. For the items (b) and (c) we refer to \cite{BDK}.
\end{proof}

For any even function  $\varphi \in \mathscr{S}(\mathbb{R})$, $\alpha>0$ and $f\in L^2(X)$ we define
$$
\varphi^*_{L,\alpha}(f)(x)=\sup_{t>0}\sup_{d(x,y)<\alpha t}|\varphi(t\sqrt{L})f(y)|,
$$ 
and 
$$
\varphi^+_{L,\alpha}(f)(x)=\sup_{t>0}|\varphi(t\sqrt{L})f(x)|.
$$
As usual, we drop the index $\alpha$ when $\alpha=1$.

The following results are taken from Proposition 2.3 and Theorem 3.1 in  \cite{SY}, respectively.
\begin{prop}
	\label{prop1-maximal}
	Let $p\in (0,1]$. Let $\varphi_1, \varphi_2\in \mathscr{\mathbb{R}}$ be even functions with $\varphi_1(0)=1$ and $\varphi_2(0)=0$ and $\alpha_1, \alpha_2>0$. Then for every $f\in L^2(X)$ we have
	\begin{equation}
	\label{eq1-prop1}
	\|(\varphi_2)^*_{L,\alpha_2}f\|_{L^p(X)}\lesi \|(\varphi_1)^*_{L,\alpha_1}f\|_{L^p(X)}.
	\end{equation}
	As a consequence, for every even function $\varphi$ with $\varphi(0)=1$ and $\alpha>0$ we have
	\begin{equation}
	\label{eq2-prop1}
	\|\varphi^*_{L,\alpha}f\|_{L^p(X)}\sim \|f^*_{L}\|_{L^p(X)}.
	\end{equation}
\end{prop}
\begin{prop}
	\label{prop2-maximal}
	Let $p\in(0,1]$. Let $\varphi\in\mathscr{S}(\mathbb{R})$ be an even function with $\varphi(0)=1$. Then we have
	\begin{equation*}
	\Big\|\varphi^*_{L}f\Big\|_{L^p(X)}\lesi \|\varphi^+_{L}f\|_{L^p(X)},
	\end{equation*}
\end{prop}

\section{Maximal function characterizations for Hardy spaces $H^p_L(X)$}\label{sect-mainproofs}
In this section we give the proofs of Theorems \ref{mainthm1} and \ref{mainthm2}.
\subsection{Proof of Theorem \ref{mainthm1}}
	By Proposition \ref{prop2-maximal} we have $H^{p}_{L, {\rm max}}(X)\equiv H^{p}_{L, {\rm rad}}(X)$ for $0<p\leq 1$ so it suffices to prove $H^p_{L,at,M}\equiv H^p_{L, \max}$ for sufficiently large $M$. The direction $H^{p}_{L,at,M}(X)\subset H^{p}_{L, {\rm max}}(X)$ follows by a similar argument to Step I in the proof of Theorem 3.5 of \cite{DHMMY} provided $p\in (0,1]$ and $M>\f{n}{2}(\f{1}{p}-1)$. The rest of the proof is devoted to the remaining direction $H^p_{\max,L}\subset H_{L,at,M}^{p}(X)$.

	Fix $f\in H^p_{L,\max}\cap L^2(X)$. We shall show that $f$ has a $(p,M)_L$-representation $\sum_j \lambda_ja_j$ with $(\sum_j|\lambda|_j^p)^{1/p} \lesssim \| f\|_{H^p_{L,\max}}$. 
	
	Let $\Phi$ be a function from Lemma \ref{lem:finite propagation}. For $M\in \mathbb{N}, M>\f{n}{2}(\f{1}{p}-1)$ we have
	\begin{equation}
	\label{Calderon forula}
	\begin{aligned}
	f=c_{\Phi,M}\int_0^\vc (t^2L)^{M}\Phi(t\sqrt{L})\Phi(t\sqrt{L})f\f{dt}{t}
	\end{aligned}
	\end{equation}
	in $L^2(X)$, where $\displaystyle c_{\Phi,M}= \Big[\int_0^\vc x^{2M}\Phi(x)\f{dx}{x}\Big]^{-1}$.	
	
	We now set
	\[
	\psi(x)=c_{\Phi,M}\int_1^\vc (tx)^{2M}\Phi^2(tx)\f{dt}{t}=c_{\Phi,M}\int_x^\vc t^{2M}\Phi^2(t)\f{dt}{t}.
	\]
	Then  $\psi \in \mathscr{S}(\mathbb{R})$ and is an even function with $\psi(0)=1$. Moreover, for $s>0$,
	\begin{equation*}
	\psi(sx)=c_{\Phi,M}\int_s^\vc (tx)^{2M}\Phi^2(tx)\f{dt}{t}.
	\end{equation*}
	This implies, for $s>0$,
	\begin{equation}
	\label{eq-psi}
	\psi(s\sqrt{L})f=c_{\Phi,M}\int_s^\vc (t\sqrt{L})^{2M}\Phi(t\sqrt{L})\Phi(t\sqrt{L})f\f{dt}{t}
	\end{equation}
	in $L^2(X)$.
	
	Setting $R_0={\rm diam}\, X/2$, we then decompose $f$ as follow
	\begin{equation}
	\label{eq-f1 f2}
	\begin{aligned}
	f&=c_{\Phi,M}\int_0^{R_0} (t^2L)^{M}\Phi(t\sqrt{L})\Phi(t\sqrt{L})f\f{dt}{t}+c_{\Phi,M}\int^\vc_{R_0} (t^2L)^{M}\Phi(t\sqrt{L})\Phi(t\sqrt{L})f\f{dt}{t}\\
	&=c_{\Phi,M}\int_0^{R_0} (t^2L)^{M}\Phi(t\sqrt{L})\Phi(t\sqrt{L})f\f{dt}{t}+\psi(R_0\sqrt{L})f\\
	&=: f_1+f_2
	\end{aligned}
	 \end{equation}
	 in $L^2(X)$.

	Define the maximal operator
	$$
	\mathbb{M}_{L}f(x)=\sup_{t>0}\sup_{d(x,y)<8t}\left[ |\psi(t\sqrt{L})f(y)| +|\Phi(t\sqrt{L})f(y)|\right].
	$$
	Then Proposition \ref{prop2-maximal} yields
	\begin{equation}
	\label{eq1-proff mainthm1}
	\|\mathbb{M}_{L}f\|_{L^p(X)}\lesi \|f\|_{H^p_{\max, L}(X)}.
	\end{equation}
	
	Let us address the term $f_2$ first. Since $R_0={\rm diam}\,X/2$, we have, for all $x,y\in X$,
	\[
	|\psi(R_0\sqrt{L})f(x)|\leq \sup_{d(z,y)<6R_0}|\psi(R_0\sqrt{L})f(z)|\le \mathbb{M}_{L}f(y).
	\]
	This implies
	\[
	\begin{aligned}
	\|f_2\|_{L^\vc(X)}&\le \inf_{y\in X}\mathbb{M}_{L}f(y)\le \mu(X)^{-1/p}\|\mathbb{M}_{L}f\|_{L^p(X)}
	\lesi  \mu(X)^{-1/p}\|f\|_{H^p_{\max, L}(X)}
	\end{aligned}
	\]
	 where in the last inequality we used \eqref{eq1-proff mainthm1}.
	 
	 Therefore, we can write $f_2= \lambda a$  with $|\lambda| \lesi \|f\|_{H^p_{\max, L}(X)}$ and $\|a\|_{L^\vc(X)}\le \mu(X)^{-1/p}$.
	 
	 It remains to decompose the term $f_1$ in terms of  $(p,M)_L$ atoms. To do this, for each $k\in \mathbb{Z}$ we set
	 $$
	 \Om_i:=\{x\in X: \mathbb{M}_{L}f(x)>2^i\}.
	 $$
	 Since $\mathbb{M}_{L}f$ is lower--continuous and $X$ is bounded, there exists $i_0$ so that $\Om_{i_0}=X$ and $\Om_{i_0+1}\neq X$. Without loss of generality we may assume that $i_0=0$. Then for each $t>0$ we define
	 \begin{equation}
	 \label{eq-Om it}
	 \Om_i^t=\begin{cases}
	 \Om_0, \ \ &i=0,\\
	 \{x: d(x,\Om_{i}^c)>4t\}, \ \ & i>0,
	 \end{cases}
	 \end{equation}
	 and $T_i^t=\Om_i^t\backslash \Om_{i+1}^t$.
	 
	 It is clear that $X=\bigcup_{i=0}^\vc T_i^t$ for each $t>0$. Hence, 
	 \begin{equation}\label{eq1-f1}
	 \begin{aligned}
	 f_1&=\sum_{i=0}^\infty c_{\Phi,M}\int_0^{R_0} (t^2L)^{M}\Phi(t\sqrt{L})\left[\Phi(t\sqrt{L})f\cdot \chi_{T_i^t}\right]\f{dt}{t}\\
	 &=:\sum_{i=0}^\infty f_1^i.
	 \end{aligned}
	 \end{equation}
	 
	 We now consider $f_1^0$ first. For $x\in X$ we have
	 \[
	 f_1^0(x)=c_{\Phi,M}\int_0^{R_0}\int_{T_0^t} K_{(t^2L)^{M}\Phi(t\sqrt{L})}(x,y)\Phi(t\sqrt{L})f(y)\dy\f{dt}{t}.
	 \] 
	 We now consider two cases: $x\in \Om_{1}^c$ and $x\in \Om_1$.
	 
	 \noindent{\textbf{Case 1: $x\in \Om_{1}^c$}.} In this situation, we can see that 
	 $$\supp K_{(t^2L)^{M}\Phi(t\sqrt{L})}(x,\cdot)\subset \{z: d(x,z)<t\}\subset T_0^t, \ \ \text{for all $t>0$}.
	 $$ 
	 Therefore,
	 \begin{equation}
	 \label{eq1-f 10}
	  \begin{aligned}
	 |f_1^0(x)|&=c_{\Phi,M}\int_0^{R_0}\int_{X} K_{(t^2L)^{M}\Phi(t\sqrt{L})}(x,y)\Phi(t\sqrt{L})f(y)\dy\f{dt}{t}\\
	 &=c_{\Phi,M}\int_0^{R_0} (t^2L)^{M}\Phi(t\sqrt{L})\Phi(t\sqrt{L})f(x)\f{dt}{t}\\
	 &=\lim_{\epsilon\to 0}c_{\Phi,M}\int_\epsilon^{R_0} (t^2L)^{M}\Phi(t\sqrt{L})\Phi(t\sqrt{L})f(x)\f{dt}{t}\\
	 &=\lim_{\epsilon\to 0}\psi(\epsilon\sqrt{L})f(x)-\psi(R_0\sqrt{L})f(x)
	 \end{aligned}
	 \end{equation}
	 where in the last inequality we used \eqref{eq-psi}.

	 On the other hand, since $x\in \Om_1^c$ we have $|\psi(s\sqrt{L})f(x)|\le 2$. This and \eqref{eq1-f 10} imply
	 \begin{equation}
	 \label{eq1-f01 outside Om1}
	 |f_1^0(x)|\le 4, \ \ \ \forall x\in \Om_1^c.
	 \end{equation}
	 
	 \noindent{\textbf{Case 2: $x\in \Om_{1}$}.} To deal with this case, we write
	 \[
	 \begin{aligned}
	 f_1^0(x)&=c_{\Phi,M}\int_0^{R_0}\int_{T_0^t} K_{(t^2L)^{M}\Phi(t\sqrt{L})}(x,y)\Phi(t\sqrt{L})f(y)\dy\f{dt}{t}\\
	 &=\int_0^{d(x,\Om_1^c)/5}\ldots+\int_{d(x,\Om_1^c)/5}^{d(x,\Om_1^c)/3}\ldots +\int_{d(x,\Om_1^c)}^{R_0}\ldots\\
	 &=:E_{1}(x)+E_2(x)+E_3(x).
	 \end{aligned}
	  \]
	 For $t\in (0,d(x,\Om_1^c)/5)$ and $y\in T_0^t$ we have $d(x,y)\ge t$. This, along with Lemma \ref{lem:finite propagation}, yields $K_{(t^2L)^{M}\Phi(t\sqrt{L})}(x,y)=0$, and hence $E_1(x)=0$.
	 
	 For the second term, using Lemma \ref{lem:finite propagation} again we have
	 \[
	 |E_2(x)|\lesi \int_{d(x,\Om_1^c)/3}^{d(x,\Om_1^c)/3}\sup_{y\in T^t_0}|\Phi(t\sqrt{L})f(y)|\f{dt}{t}.
	 \]
	 From the definition of the set $T^t_0$ it is easy to see that for each $y\in T^t_0$ we can find $z\in \Om_1^c$ so that $d(y,z)<6t$. For each such $z$ we have, since $z\in \Om_1^c$,
	 \[
	 |\Phi(t\sqrt{L})f(y)|\le \mathbb{M}_Lf(z)\le 2.
	 \]

	 Therefore we obtain $|E_2(x)|\lesssim 1$.

	 For the last term $E_3(x)$, we observe that for $t>d(x,\Om_1^c)/3$ and $x\in \Om_1$ we have
	 $$
	 \supp K_{(t^2L)^{M}\Phi(t\sqrt{L})}(x,\cdot)\subset \{z: d(x,z)<t\}\subset T_0^t, \ \ \text{for all $t>0$}.
     $$	  
	 Hence, arguing similarly to \eqref{eq1-f 10} we come up with
	 \[
	 E_3(x)=\psi(s_1\sqrt{L})f(x)-\psi(R_0\sqrt{L})f(x)
	 \]
	 where $s_1=d(x,\Om_1^c)/3$.
	 
	 Note that for $s>d(x,\Om_1^c)/3$ we can find $z\in \Om_1^c$ so that $d(x,z)<3s$. Hence,
	 \[
	 |\psi(t\sqrt{L})f(x)|\le \mathbb{M}_Lf(z)\le 2,
	 \]
	 since $z\in \Om_1^c$.
	 
	As a consequence, we have $|E_3(x)|\lesi 1$. We now take all estimates $E_1(x), E_2(x)$ and $E_3(x)$ into account to find that
	\begin{align}\label{thm1eq1}
	|f_1^0(x)|\lesi 1, \ \ \ \forall x\in \Om_1.
	\end{align}
	This, along with \eqref{eq1-f01 outside Om1}, implies that \eqref{thm1eq1} in fact holds for every $x\in X$.
	Then we have
	\[
	\begin{aligned}
	|f_1^0(x)|&\lesi \mu(X)^{-1/p}\mu(\Om_0)^{1/p}\\
	&\lesi \mu(X)^{-1/p}\sum_{i=0}^\vc 2^i\mu(\Om_i)^{1/p}\sim \mu(X)^{-1/p}\|\mathbb{M}_Lf\|_{L^p(X)}\\
	&\lesi \mu(X)^{-1/p}\|f\|_{H^p_{L,\max}(X)}.	
	\end{aligned}
		\]
Hence, we can write $f_1^0= \lambda_1^0 a_1^0$ so that $|\lambda_1^0| \lesi \|f\|_{H^p_{\max, L}(X)}$ and $\|a_1^0\|_{L^\vc(X)}\le \mu(X)^{-1/p}$.	 
	 
We now take care of the term $f_1^i$ with $i>0$. To do this, for each $i>0$ we apply a covering lemma in \cite{CW} (see also \cite[Lemma 5.5]{DKP}) to obtain a collection of balls $\{B_{i,k}:=B(x_{B_{i,k}},r_{B_{i,k}}): x_{B_{i,k}}\in \Om_i, r_{B_{i,k}}=d(x_{B_{i,k}},\Om_i^c)/2, k=1,\ldots\}$ so that
	\begin{enumerate}[{\rm (i)}]
		\item $\displaystyle \Om_i=\cup_k B(x_{B_{i,k}},r_{B_{i,k}})$;
		\item $\displaystyle  \{B(x_{B_{i,k}},r_{B_{i,k}}/5)\}_{k=1}^\vc$ are disjoint.
	\end{enumerate}
For each $i, k\in \mathbb{N}^+$  and $t>0$ we set $B^t_{i,k}=B(x_{i,k},r_{B_{i,k}}+2t)$ which is a ball having the same center as $B_{i,k}$ with radius being $2t$ greater than the radius of  $B_{i,k}$. Then, for each $i, k\in \mathbb{N}^+$  and $t>0$, we set
\[
R_{i,k}^t=\begin{cases}
T_i^t\cap B^t_{i,k}, \ \ &\text{if} \ \ T_i^t\cap B_{i,k}\ne \emptyset\\
0,  \ \ &\text{if} \ \ T_i^t\cap B_{i,k}= \emptyset,
\end{cases}
\]
and 
\begin{equation}
\label{eq-Eit}
E_{i,k}^t=R_{i,k}^t\backslash \cup_{\ell>k}R_{i,k}^t.
\end{equation}
It is easy to see that for each $i\in \mathbb{N}^+$ and $t>0$ we have
\[
T_i^t =\bigcup_{k\in \mathbb{N}^+}E_{i,k}^t.
\]
Hence, from \eqref{eq1-f1} we have, for $i\in \mathbb{N}^+$,
	\[
	\begin{aligned}
	f^i_1&=\sum_{k\in \mathbb{N}^+} c_{\Phi,M}\int_0^{R_0} (t^2L)^{M}\Phi(t\sqrt{L})\left[\Phi(t\sqrt{L})f\cdot \chi_{E_{i,k}^t}\right]\f{dt}{t}
	\end{aligned}
	\]
	and set $a_{i,k}=0$ if $E_{i,k}^t=\emptyset$.
	
	We now define $\lambda_{i,k}=2^i\mu(B_{i,k})^{1/p}$ and $a_{i,k}=L^Mb_{i,k}$ where
	\begin{equation}\label{eq-bik}
	b_{i,k}=\f{c_{\Phi,M}}{\lambda_{i,k}}\int_0^{R_0} t^{2M}\Phi(t\sqrt{L})\left[\Phi(t\sqrt{L})f\cdot \chi_{E_{i,k}^t}\right]\f{dt}{t}.
	\end{equation}
	Then it can be seen that 
	$$
	\begin{aligned}
	f_1&=\sum_{i\in \mathbb{N}^+}f_1^i=\sum_{i,k\in \mathbb{N}^+}\lambda_{i,k}a_{i,k}
	\end{aligned}
	$$
	in $L^2(X)$; moreover,
	\[
	\begin{aligned}
	\sum_{i,k\in \mathbb{N}^+}|\lambda_{i,k}|^p&=\sum_{i,k\in \mathbb{N}^+} 2^{ip}\mu(B_{i,k})\lesi \sum_{i\in \mathbb{N}^+} 2^{ip}\mu(\Om_i)
	\lesi \|\mathbb{M}_Lf\|^p_{L^p(X)}\lesi \|f\|^p_{H^p_{L,\max}(X)}.
	\end{aligned}
	\]
	Therefore, it suffices to prove that each $a_{i,k}\ne 0$ is a $(p,M)_L$ atom associated to the ball $B^*_{i,k}:=8B_{i,k}$. Indeed, if $r_{B_{i,k}}<t/2$, then we have $d(x_{B_{i,k}}, \Om_i^c)=2r_{B_{i,k}}<t$. Therefore, 
	$$
	B^t_{i,k}=B(x_{B_{i,k}}, r_{B_{i,k}}+2t) \subset \{x: d(x,\Om_i^c)<4t\}.
	$$
	This implies that $R_{i,k}^t:=T_i^t\cap B^t_{i,k}=\emptyset$. Hence, if $a_{i,k}\ne 0$, then $r_{B_{i,k}}\ge t/2$. This, along with \eqref{eq-bik} and Lemma \ref{lem:finite propagation}, implies that 
	\[
	\supp L^mb_{i,k}\subset B^*_{i,k}, \ \ \ \forall m=0,1,\ldots, M.
	\]
	It remains to show that
	\[
	|L^mb_{i,k}|_{L^\vc(X)}\lesi r_{B_{i,k}}^{2(M-m)}\mu(B_{i,k})^{-1/p}, \ \ \ \forall m=0,1,\ldots, M.
	\]
	For $m=0,1,\ldots, M-1$, since $r_{B_{i,k}}\ge t/2$ as $a_{i,k}\ne 0$, we have
	\[
	|L^mb_{i,k}(x)|\le \f{c_{\Phi,M}}{\lambda_{i,k}}\int_0^{2r_{B_{i,k}}} t^{2(M-m)}\Big|(t^2L)^m\Phi(t\sqrt{L})\left[\Phi(t\sqrt{L})f\cdot \chi_{E_{i,k}^t}\right](x)\Big|\f{dt}{t}.
	\]
	This along with Lemma \ref{lem:finite propagation} implies that
	\[
	\begin{aligned}
	|L^mb_{i,k}(x)|&\lesi \f{1}{\lambda_{i,k}}\int_0^{2r_{B_{i,k}}} t^{2(M-m)}\int_{E_{i,k}^t} |\Phi(t\sqrt{L})f(y)|\dy\f{dt}{t}\\
	&\lesi \f{1}{\lambda_{i,k}}\int_0^{2r_{B_{i,k}}} t^{2(M-m)}\int_{T_i^t} |\Phi(t\sqrt{L})f(y)|\dy\f{dt}{t}.
	\end{aligned}
		\]
		Note that for each $y\in T_i^t$ there exists $z\in \Om_{i+1}^c$ so that $d(y,z)<4t$. Hence,
		\[
		|\Phi(t\sqrt{L})f(y)|\le \mathbb{M}_Lf(z)\le 2^{i+1}, \ \ \ \forall y\in T_i^t.
		\]
	Therefore, for all $m=0,1,\ldots, M-1$,
	\[
	\begin{aligned}
	|L^mb_{i,k}(x)|	&\lesi \f{2^i}{\lambda_{i,k}}\int_0^{2r_{B_{i,k}}} t^{2(M-m)}\f{dt}{t}\lesi \mu(B_{i,k})^{-1/p}.
	\end{aligned}
	\]	
	For $m=M$, we have
	\[
	L^Mb_{i,k}(x)=\f{c_{\Phi,M}}{\lambda_{i,k}}\int_0^{2r_{B_{i,k}}} (t^2L)^M\Phi(t\sqrt{L})\left[\Phi(t\sqrt{L})f\cdot \chi_{E_{i,k}^t}\right](x)\f{dt}{t}.
	\]
	From \eqref{eq-Eit} we have
	\[
	E_{i,k}^t=(T_i^t\cap B_{i,k}^t)\backslash (T_i^t\cap F_{i,k}^t)
	\]
	where $F_{i,k}^t:=\cup_{\ell>k}B_{i,\ell}^t=\{x: d(x, \cup_{\ell>k}B_{i,\ell})<2t\}$.
	
	Hence, we can rewrite
	\begin{equation}
	\label{eq-LMb}
	\begin{aligned}
	L^Mb_{i,k}(x)=&\,\f{c_{\Phi,M}}{\lambda_{i,k}}\int_0^{2r_{B_{i,k}}}\int_{T_i^t\cap B_{i,k}^t} K_{(t^2L)^M\Phi(t\sqrt{L})}(x,y)\Phi(t\sqrt{L})f(y)\dy\f{dt}{t}\\
	& - 	\f{c_{\Phi,M}}{\lambda_{i,k}}\int_0^{2r_{B_{i,k}}}\int_{T_i^t\cap F^t_{i,k}\cap B_{i,k}^t} K_{(t^2L)^M\Phi(t\sqrt{L})}(x,y)\Phi(t\sqrt{L})f(y)\dy\f{dt}{t}\\
	=&\, \f{c_{\Phi,M}}{\lambda_{i,k}}I_{i,k}(x)+\f{c_{\Phi,M}}{\lambda_{i,k}}J_{i,k}(x).
	\end{aligned}
	\end{equation}
	We have the following result whose proof will be given after the proof of the theorem.
	\begin{lem}\label{lem-mainproof}
		There exists a constant $C>0$ so that for all $i, k\in \mathbb{N}^+$ and $x\in X$ we have
		\begin{equation}
		\label{eq-estimate I, J}
		|I_{i,k}(x)|+|J_{i,k}(x)|\lesi 2^i,
		\end{equation}
	where $I_{i,k}$ and $J_{i,k}$ have been defined in \eqref{eq-LMb}.
	\end{lem} 
	We now just substitute \eqref{eq-estimate I, J} into \eqref{eq-LMb} to conclude that 
	$$|L^Mb_{i,k}(x)|\lesi \mu(B_{i,k})^{-1/p}.
	$$
	which concludes the proof of Theorem \ref{mainthm1}.


\bigskip

We now give the proof for Lemma \ref{lem-mainproof}.
\begin{proof}
	[Proof of Lemma \ref{lem-mainproof}:] 
	For any subset $U$ of $X$ and for each $t>0$ we define $U^t:=\{x: d(x,U)<2t\}$.
	Now let  $U$ and $V$ be any two subsets of $X$. For each $s\in (0,R_0]$ and $i\in \mathbb{N}^+$ we define
	\[
	g_s(x)=c_{\Phi,M}\int_0^{s} (t^2L)^{M}\Phi(t\sqrt{L})\left[\Phi(t\sqrt{L})f\cdot \chi_{T_i^t}\right]\f{dt}{t}, \ \  x\in X.
	\]
   \[
   g_{U,s}(x)=c_{\Phi,M}\int_0^{s} (t^2L)^{M}\Phi(t\sqrt{L})\left[\Phi(t\sqrt{L})f\cdot \chi_{T_i^t\cap U^t}\right]\f{dt}{t}, \ \  x\in X.
   \]
       
    and
     \[
    g_{U,V,s}(x)=c_{\Phi,M}\int_0^{s} (t^2L)^{M}\Phi(t\sqrt{L})\left[\Phi(t\sqrt{L})f\cdot \chi_{T_i^t\cap U^t\cap V^t}\right]\f{dt}{t}, \ \  x\in X.
    \]
    
    We claim that \eqref{eq-estimate I, J} is a consequence of the following three estimates.
	\begin{equation}
	\label{eq-gs}
	|g_s(x)|\lesi 2^i, 
	\end{equation}
   \begin{equation}
   \label{eq-gMs}
   |g_{U,s}(x)|\lesi 2^i,
   \end{equation}
  \begin{equation}
    \label{eq-gMNs}
    |g_{U,V,s}(x)|\lesi 2^i,
    \end{equation}
    for any $U,V\subset X$ and  $s\in (0,R_0], i\in \mathbb{N}^+$ and $x\in X$
    
Indeed by firstly applying \eqref{eq-gMs} for $U=B_{i,k}$ we obtain $|I_{i,k}(x)|\lesi 2^i$ for all $x\in X$. Secondly by applying \eqref{eq-gMNs} for  $U=B_{i,k}$  and  $V=F_{i,k}$ we get $|J_{i,k}(x)|\lesi 2^i$. Thus \eqref{eq-estimate I, J} holds. 
	
It remains to show \eqref{eq-gs}--\eqref{eq-gMs}.

	We begin with \eqref{eq-gs}.
	Indeed, we now consider two cases: $x\in \Om_{i+1}^c$ and $x\in \Om_{i+1}$.
	\medskip
	
	\noindent\textbf{Case 1: $x\in \Om_{i+1}^c$.} There are two subcases $s>d(x,\Om_i^c)/3$ and $s\le d(x,\Om_i^c)/3$. We just consider the first case, since the latter is similar and even easier.
	
	We write
	\[
	\begin{aligned}
	g_s(x)=& \, c_{\Phi,M}\int_0^{d(x,\Om_i^c)/5}\int_{T_i^t} K_{(t^2L)^{M}\Phi(t\sqrt{L})}(x,y)\Phi(t\sqrt{L})f(y)\dy\f{dt}{t}\\
	&+c_{\Phi,M}\int_{d(x,\Om_i^c)/5}^{d(x,\Om_i^c)/3} \int_{T_i^t} K_{(t^2L)^{M}\Phi(t\sqrt{L})}(x,y)\Phi(t\sqrt{L})f(y)\dy\f{dt}{t}\\
	&+c_{\Phi,M}\int_{d(x,\Om_i^c)/3}^s \int_{T_i^t} K_{(t^2L)^{M}\Phi(t\sqrt{L})}(x,y)\Phi(t\sqrt{L})f(y)\dy\f{dt}{t}\\
	=:&\, A_1(x)+A_2(x)+A_3(x).
	\end{aligned}
	\]
	For the first term, we can see that $B(x,t)\subset T_i^t$ as $t\in (0, d(x,\Om_i^c)/5)$. Hence, by Lemma \ref{lem:finite propagation} we find that
	\[
	\begin{aligned}
	A_1(x)&=c_{\Phi,M}\int_0^{d(x,\Om_i^c)/5}\int_{X} K_{(t^2L)^{M}\Phi(t\sqrt{L})}(x,y)\Phi(t\sqrt{L})f(y)\dy\f{dt}{t}\\
	&=c_{\Phi,M}\int_0^{d(x,\Om_i^c)/5}(t^2L)^{M}\Phi(t\sqrt{L})\Phi(t\sqrt{L})f(x)\f{dt}{t}\\
		 &=\lim_{\epsilon\to 0}\psi(\epsilon\sqrt{L})f(x)-\psi(s_2\sqrt{L})f(x)
	\end{aligned}
	\]
	where $s_2=d(x,\Om_i^c)/5$.
	
	Note that $x\in \Om_{i+1}^c$. Hence, $|\psi(t\sqrt{L})f(x)|\le 2^{i+1}$ for all $t>0$. As a consequence, we have $A_1(x)\le 2^{i+1}$. 

    For the second term $A_2(x)$, using Lemma \ref{lem:finite propagation} we obtain
    \[
    |A_2(x)|\lesi \int_{d(x,\Om_i^c)/5}^{d(x,\Om_i^c)/3} \sup_{y\in T_i^t}|\Phi(t\sqrt{L})f(y)|\f{dt}{t}.
    \]
    Moreover, since $B(y,8t)\cap \Om_{i+1}^c\ne \emptyset$ as $y\in T_i^t$ and $d(x,\Om_i^c)/5\le t \le d(x,\Om_i^c)/3$, we have
    \[
    \sup_{y\in T_i^t}|\Phi(t\sqrt{L})f(y)|\le 2^{i+1}, \ \ \forall d(x,\Om_i^c)/5\le t \le d(x,\Om_i^c)/3.
    \]
    Hence, $A_2(x)\lesi 2^i$.
    
    The last term $A_3(x)$ is zero, since in this situation we have $B(x,t)\cap T_i^t=\emptyset$. Gathering all estimates of $A_1(x), A_2(x)$ and $A_3(x)$ we arrive at $|g_s(x)|\lesi 2^i$.
    
    \medskip
    
    \noindent\textbf{Case 2: $x\in \Om_{i+1}$.} We only consider the case $s>d(x,\Om_{i+1}^c)$. The case $s\le d(x,\Om_{i}^c)$ can be done similarly.
    
    In this case, firstly we consider situation $d(x,\Om_{i+1}^c)/3<d(x,\Om_{i}^c)/4$. Then we split the integral in the expression of $g_s$ as follows
    \[
    \begin{aligned}
    g_s(x)=& \,\int_0^{d(x,\Om_{i+1}^c)/5}\ldots+\int_{d(x,\Om_{i+1}^c)/5}^{d(x,\Om_{i+1}^c)/3}\ldots+\int_{d(x,\Om_{i+1}^c)/3}^{d(x,\Om_{i}^c)/4}\ldots\\
    &+\int_{d(x,\Om_{i}^c)/4}^{d(x,\Om_{i}^c)/3}\ldots+\int_{d(x,\Om_{i}^c)/3}^{s}\ldots=: \sum_{\ell=1}^5B_\ell(x).
    \end{aligned}
    \]
    Arguing as in the first case, we have $B_1(x)=B_5(x)=0$. The terms $B_2(x)$ and $B_4(x)$ can be dealt with in a similar way to  $A_2(x)$ so that
    \[
    |B_2(x)|+|B_4(x)|\lesi 2^i.
    \]
    For the term $B_3(x)$, we note that by Lemma \ref{lem:finite propagation} we have,
    \[
    \supp K_{(t^2L)^{M}\Phi(t\sqrt{L})}(x,y)(x,\cdot)\subset B(x,t)\subset T_i^t, \ \ \ d(x,\Om_{i+1}^c)/3<t<d(x,\Om_{i}^c)/4.
    \]
    At this stage, arguing similarly to the estimate of $A_1(x)$ in Case 1, we find that
    \[
    |B_3(x)|\lesi 2^i.
    \]
    If $d(x,\Om_{i+1}^c)/3\ge d(x,\Om_{i}^c)/4$, we split the integral in the expression of $g_s$ as follows:
    \[
    \begin{aligned}
    g_s(x)=& \,\int_0^{d(x,\Om_{i+1}^c)/5}\ldots
    +\int_{d(x,\Om_{i+1}^c)/5}^{d(x,\Om_{i}^c)/4}\ldots+\int_{d(x,\Om_{i}^c)/4}^{d(x,\Om_{i}^c)/3}\ldots+\int_{d(x,\Om_{i}^c)/3}^{s}\ldots.
    \end{aligned}
    \]
   Then we use the argument as above to dominate $|g_s(x)|$ by a multiple of  $2^i$. 
   Hence, this completes the proof of \eqref{eq-gs}.

We turn to the proof of \eqref{eq-gMs}.
   
   Now if $x\in U$, then 
   \[
   \supp K_{(t^2L)^{M}\Phi(t\sqrt{L})}(x,y)(x,\cdot)\subset B(x,t)\subset U^t, \ \ \forall t>0.
   \]
   Hence, $g_{U,s}(x)=g_s(x)$, and by \eqref{eq-gs} we have $|g_{U,s}(x)|\lesi 2^i$. 
   
   Otherwise, if $x\notin U$, there two cases $s>d(x,U^c)$ and $s\le d(x,U^c)$. We will only consider the first case $s>d(x,U^c)$,  since the second case can be done similarly.
   Assuming $s>d(x,U^c)$, we now break $g_{U,s}$ into 3 terms
   \[
   \begin{aligned}
   g_{U,s}(x)=&\,c_{\Phi,M}\int_0^{d(x,U^c)/3}\int_{T_i^t\cap U^t}K_{ (t^2L)^{M}\Phi(t\sqrt{L})}(x,y)\Phi(t\sqrt{L})f(y)\dy\f{dt}{t}\\
   &+   c_{\Phi,M}\int_{d(x,U^c)/3}^{d(x,U^c)}\int_{T_i^t\cap U^t}K_{ (t^2L)^{M}\Phi(t\sqrt{L})}(x,y)\Phi(t\sqrt{L})f(y)\dy\f{dt}{t}\\
   &+   c_{\Phi,M}\int_{d(x,U^c)}^s\int_{T_i^t\cap U^t}K_{ (t^2L)^{M}\Phi(t\sqrt{L})}(x,y)\Phi(t\sqrt{L})f(y)\dy\f{dt}{t}\\
   =&\, D_1(x)+D_2(x)+D_3(x).
   \end{aligned}
   \]
   The first term $D_1(x)=0$, since 
   \[
   \supp K_{(t^2L)^{M}\Phi(t\sqrt{L})}(x,y)(x,\cdot)\subset B(x,t)\subset (U^t)^c.
   \]
   The second term can be estimated similarly to the term $A_2(x)$ so that $|D_2(x)|\lesi 2^i$. 
   
   For the last term, we note that as $t>d(x,U^c)$ we have
   \[
   \supp K_{(t^2L)^{M}\Phi(t\sqrt{L})}(x,y)(x,\cdot)\subset B(x,t)\subset U^t.
   \]
    Hence,
    \[
    D_3(x)=c_{\Phi,M}\int_{d(x,U^c)}^s\int_{T_i^t}K_{ (t^2L)^{M}\Phi(t\sqrt{L})}(x,y)\Phi(t\sqrt{L})f(y)\dy\f{dt}{t}=g_s(x)-g(s_3)(x)
    \]
    where $s_3=d(x,U^c)\le R_0$.
    
    Using \eqref{eq-gs} we obtain $|D_3(x)|\lesi 2^i$ and this proves \eqref{eq-gMs}.
    
   The estimate  \eqref{eq-gMNs} can be done by repeating the argument for \eqref{eq-gMs}. This completes the proof of Lemma \ref{lem-mainproof}. 
   \end{proof}    

\begin{rem}
	The argument in the proof of Theorem \ref{mainthm1} also shows that the result in Theorem \ref{mainthm1} is still true if we define  the \textbf{non-tangential}  maximal function and  the \textbf{radial} maximal function by
	\[
	f^*_{L}(x)=\sup_{0<t<{\rm diam^2\,}X}\sup_{d(x,y)<t}|e^{-tL}f(y)|
	\]
	and
	\[
	f^+_{L}(x)=\sup_{0<t<{\rm diam^2\,}X}|e^{-tL}f(x)|,
	\]
	respectively.
\end{rem}
\subsection{Proof of Theorem \ref{mainthm2}}
	
	Due to the validity of Theorem \ref{mainthm1}, we need only to show that $H^p_{CW}(X)\equiv H^p_{L,\rad}(X)$ for $\f{n}{n+\delta}<p\le 1$. Now the inclusion $H^p_{CW}(X)\hookrightarrow H^p_{L,\rad}(X)$ is standard and will be omitted. For the reverse inclusion, since $H^p_{L,\rad}(X)\equiv H^p_{L,at,M}(X)$, again from Theorem \ref{mainthm1}, it suffices to prove $H^p_{L,at, M}(X)\hookrightarrow H^p_{CW}(X)$. 
	
	Now if $a$  is a $(p,M)_L$-atom associated to a ball $B$, then there are two cases. 
	If $a=\mu(X)^{-1/p}$ then $a$ is also a $p$-atom and hence $a\in H^p_{CW}(X)$.
	Otherwise, we have $a=L^Mb$ and by using (A4)  and arguing similarly to the proof of Lemma 9.1 in \cite{HLMMY} we have
	\[
	\int a(x)\dx =0.
	\]
	Moreover parts (ii) and (iii) of Definition \ref{def: L-atom}  shows that $\supp a \subset B$ and $\|a\|_{L^\vc}\le \mu(B)^{-1/p}$. It follows that $a$ is also a $p$-atom and hence $a\in H^p_{CW}(X)$. As a consequence, $H^p_{L,at, M}(X)\hookrightarrow H^p_{CW}(X)$ and this completes our proof of Theorem \ref{mainthm2}.
 \section{Applications}\label{sect-app}
 
 In this section we apply Theorems \ref{mainthm1} and \ref{mainthm2} to give maximal characterizations of atomic Hardy spaces for various differential operators on domains. 

\subsection{Second-order elliptic operators with Neumann boundary conditions}\label{sect Neuman}
Let $A: \Rn \to M_n(\mathbb{R})$ be a real matrix-valued function and define
\[
\|A\|_\vc =\sup_{x\in \Rn, |\xi|=|\eta|=1}|A(x)\xi\cdot \eta|.
\]
We assume that $A$ is symmetric and satisfies the following conditions for all $x\in \Rn, \xi\in \mathbb{R}^n$:
\begin{equation}\label{eq-condition A}
\|A\|_\vc\le \Lambda^{-1} \ \ \text{and} \ \ \,A(x)\xi\cdot\xi \ge \Lambda |\xi|^2
\end{equation}
for some $\Lambda>0$.

Let $\Om$ be an connected open bounded/unbounded domain in $\mathbb{R}^n$ satisfying a doubling property. We 
do not assume any smoothness assumption on the boundary of $\Om$ unless it is implied by other assumptions. Denote by $L_N$ the maximal-accretive operator on $L^2(\Om)$ with largest domain $D(L_N)\subset W^{1,2}(\Omega)$ such that
\begin{equation}\label{eq-formula L}
\langle L_Nf, g \rangle =\int_{\Omega} A\nabla f\cdot \nabla g, \ \ \ \forall f\in D(L_N), \ \ g\in W^{1,2}(\Om).
\end{equation}
We then have the following.
\begin{thm}\label{mainth-Neumann}
	Assume that the kernel $p_t(x,y)$ of $e^{-tL_N}$ satisfies the following conditions:
	\begin{enumerate}
		\item[(N1)] There exists $C, C>0$ so that
		\[
		|p_t(x,y)|\le \f{C}{\mu(B_\Omega(x,\sqrt{t}))}\exp\Big(-\f{|x-y|^2}{ct}\Big)
		\]
		for all $0<t<{\rm diam}\,X$ and $x,y\in \Omega$, where $B_\Omega(x,r)=B(x,r)\cap \Om$. 
		\item[(N2)] There exist $\gamma\in (0,1]$ and $C,c>0$ so that
		\[
		|p_t(x,y)-p_t(x,y')|\le \Big(\f{|y-y'|}{\sqrt{t}}\Big)^\gamma\f{C}{\mu(B_\Omega(x,\sqrt{t}))}\exp\Big(-\f{|x-y|^2}{ct}\Big)
		\]
		for all $0<t<{\rm diam}\,\Om$ and  $x,y,y'\in \Om$ so that $|y-y'|<\sqrt{t}/2$.
	\end{enumerate}
Then we have
\begin{equation}
\label{eq-Neuman}
H^p_{L_N, \max}(\Om)\equiv H^p_{L_N, \rad}(\Om)\equiv H^p_{CW}(\Om), \ \ \ \f{n}{n+\gamma}<p\le 1.
\end{equation}
\end{thm}
\begin{rem}
We note that when $\Om$ is a strongly Lipschitz domain, the equivalence \eqref{eq-Neuman} was obtained in \cite{AR} for $p=1$ and $|\Om|=\vc$. Hence, in the case $|\Om|<\vc$ our result is new even for $p=1$.
\end{rem}

\begin{proof}[Proof of Theorem \ref{mainth-Neumann}]
	We need only to prove the case $|\Om|<\vc$, since the case $|\Om|=\vc$ is similar and easier.
	Now it is well-known that 
	for every $x\in \Omega$ and $t>0$, we have
	\begin{align*}
	\displaystyle \int_\Omega p_t(x,y)\dy =1.
	\end{align*}
	Therefore, $L_N$ satisfies (A1)-(A4) and we may invoke Theorem \ref{mainthm2} to conclude our proof.
\end{proof}

\subsection{Second-order elliptic operators with Dirichlet boundary conditions}
Let $A$ and $\Om$ be as in Subsection \ref{sect Neuman}. 
Denote by $L_D$ the maximal-accretive operator on $L^2(\Om)$ with largest domain $D(L_D)\subset W^{1,2}_0(\Omega)$ such that
\begin{equation}\label{eq-formula L}
\langle L_Df, g \rangle =\int_{\Omega} A\nabla f\cdot \nabla g, \ \ \ \forall f\in D(L_D), \ \ g\in W^{1,2}_0(\Om).
\end{equation}

We shall consider the atomic spaces defined by Miyachi \cite{Mi}.
\begin{defn}[Hardy spaces of Miyachi]
Let $p\in (0,1]$. A bounded, measurable function $a: \Om\to \mathbb{R}$ is called an $H^p_{Mi}(\Om)$-atom if 
\begin{enumerate}[\upshape (i)]
	\item $a$ is supported in a ball $B\subset \Om$;
	\item $\|a\|_{L^\vc(\Om)}\le |B|^{-1/p}$;
	\item either $2B\subset \Om$ and  $4B\cap \partial\Om\ne \emptyset$, or $4B\subset \Om$ and 
	$$\displaystyle \int x^\alpha a(x)dx=0$$ 
	for all multi--indices $\alpha$ with $|\alpha|\le [n(1/p-1)]$.
\end{enumerate}
The Hardy space $H^p_{Mi}(\Om)$ is defined as the set of all $f\in \mathscr{S}'$ such that
\[
f=\sum_j\lambda_j a_j
\]
where $a_j$ are $H^p_{Mi}(\Om)$-atoms and $\lambda_j$ are scalars with $\sum_{j}|\lambda_j|^p<\vc$.
We also set
\[
\|f\|^p_{H^p_{Mi}(\Om)}=\inf\Big\{\sum_{j}|\lambda_j|^p: f=\sum_j\lambda_j a_j \Big\}
\]
where the infimum is taken over all such decompositions.
\end{defn}

Our main result in this section is the following:
\begin{thm}
	\label{mainth-Dirichlet}
Assume that the kernel $p_t(x,y)$ of $e^{-tL_D}$ satisfies the following conditions:
\begin{enumerate}
	\item[(D1)] There exists $C, C>0$ so that
	\[
	|p_t(x,y)|\le \f{C}{\mu(B_\Omega(x,\sqrt{t}))}\exp\Big(-\f{|x-y|^2}{ct}\Big)
	\]
	for all $t>0$ and $x,y\in \Omega$, where $B_\Omega(x,r)=B(x,r)\cap \Om$. 
	\item[(D2)] There exist $\gamma\in (0,1]$ and $C,c>0$ so that
	\[
	|p_t(x,y)-p_t(x,y')|\le \Big(\f{|y-y'|}{\sqrt{t}}\Big)^\gamma\f{C}{\mu(B_\Omega(x,\sqrt{t}))}\exp\Big(-\f{|x-y|^2}{ct}\Big)
	\]
	for all $0<t<{\rm diam}\,\Om$, $x,y,y'\in \Om$ so that $|y-y'|<\sqrt{t}/2$.
\end{enumerate}
Then we have
 \[
 H^p_{L_D,\max}(\Om) \equiv H^p_{L_D,\rad}(\Om)\equiv H^p_{Mi}(\Om), \ \ \f{n}{n+\gamma}<p\le 1.
 \]
\end{thm}

\begin{rem}
Some comments are in order.
\begin{enumerate}[(a)]
\item In the particular case when 
$\Om$ is $\Rn$ or Lipschitz domains, the conditions (D1) and (D2) are always satisfied. See  \cite{AR}.  
	\item Let $\phi \in C^\vc_c(B(0,1))$  be a non--negative radial function such that $\int \phi(x)dx=1$. It was proved in \cite{Mi} that the Hardy spaces $H^p_{Mi}(\Om)$ can be characterized in terms of maximal functions of the form
\[
f^+(x)=\max_{0<t<\delta(x)/2}|\phi_t\ast f(x)|
\]
where $\delta(x)=d(x,\Om^c)$ and $\phi_t(x)=t^{-n}\phi(x/t)$. More precisely, we have
\[
\|f\|_{H^p_{Mi}(\Om)}\sim \|f^+\|_{L^p(\Om)}, \ \ 0<p\le 1.
\]
In this sense, our results give new maximal characterizations for the Hardy spaces $H^p_{Mi}(\Om)$.

\item The Hardy space $H^p_{Mi}(\Om)$ is closely related to the Hardy space $H^p_r(\Om)$ defined by
\[
H^p_r(\Om)=\{f\in \mathscr{S}': \ \text{there exists $F\in H^p(\Rn)$ so that $F|_{\Om}=f$}\}
\]
with the norm
\[
\|f\|_{H^p_r(\Om)}=\inf\{\|F\|_{H^p_r(\Rn)}: F\in H^p(\Rn), F|_{\Om}=f\}.
\]
It is well-known that if $\Om$ is a strongly Lipschitz domain such that either $\Om$ is bounded or $\Om^c$ is unbounded (see Subsection 4.3 for the precise definition), then $H^p_{Mi}(\Om)\equiv H^p_{r}(\Om)$ for $0<p\le 1$. See for instance \cite{CKS}. For such domains, Theorem \ref{mainth-Dirichlet} implies that $H^p_{L,\rad}(\Om)\equiv H^p_{L,\max}(\Om)\equiv H^p_{r}(\Om)$. This gives a positive answer to the open question in \cite{DHMMY} (mentioned in Section 5.1 of that article).
\end{enumerate}
\end{rem}

Before coming to the proof of Theorem \ref{mainth-Dirichlet}, we need the following two technical results.
\begin{lem}\label{lem1-Dirichlet}
	Let $x_0\in \Om$ and $r>0$ so that $B(x_0,2r)\subset \Om$. Then we have
	\begin{equation}\label{eq-lem1D}
	\Big|\int_{B(x_0,2r)}q_t(x,y)dx\Big| \lesi \f{\sqrt{t}}{r}e^{-\f{r^2}{ct}}
	\end{equation}
	for all $t>0$ and all $y\in B(x_0,r)$ where $q_t(x,y)$ is the kernel of $tL_De^{-tL_D}$.
\end{lem}
\begin{proof}
	Let $\psi \in C^\vc_c(x_0,2r)$ so that $\psi=1$ on $B(x_0,3r/2)$ and $|\nabla \psi|\lesi 1/r$. Then we have
	\[
	\begin{aligned}
	\int_{B(x_0,2r)}q_t(x,y)dx&=\int_{B(x_0,2r)}q_t(x,y)\psi(x)dx
	 +\int_{B(x_0,2r)}q_t(x,y)[1-\psi(x)]dx
	=:I_1+I_2.
	\end{aligned}
	\]
The Gaussian upper bound of $q_t(x,y)$ and the support condition of $(1-\psi)$ gives
	\[
	\begin{aligned}
	|I_2|
	\lesi \exp\Big(-\f{r^2}{c't}\Big)\int_{B(x_0,2r)\backslash B(x_0,3r/2)} \f{1}{t^{n/2}}\exp\Big(-\f{|x-y|^2}{2ct}\Big)dx
	\lesi \f{\sqrt{t}}{r}e^{-\f{r^2}{c't}}.
	\end{aligned}
	\]
	
	For the term $I_1$ we first note that $q_t(x,y)=tL_D[p_t(\cdot,y)](x)$ and $p_t(\cdot,y)\in D(L_D)$. See for example \cite{Ou}. Hence, from \eqref{eq-formula L} and \eqref{eq-condition A}, coupled with the support property of $\nabla \psi$, we have
	\begin{equation}\label{eq-I2-Dirichlet}
	\begin{aligned}	
	|I_1|&=t\Big|\int_{\Om} A\nabla_x p_t(x,y)\cdot \nabla \psi(x)dx\Big|\\
	&\lesi \f{t}{r}\int_{B(x_0,2r)\backslash B(x_0,3r/2)} |\nabla_x p_t(x,y)|dx\\
	&\lesi tr^{n/2-1}\Big(\int_{B(x_0,2r)\backslash B(x_0,3r/2)} |\nabla_x p_t(x,y)|^2dx\Big)^{1/2}\\
	\end{aligned}
	\end{equation}
	Arguing similarly to \cite[Lemma 3]{DM} we find that there exists $\alpha>0$ so that
	\[
	\Big(\int_{\Om}|\nabla_x p_t(x,y)|^2e^{\f{|x-y|^2}{\alpha t}}dx\Big)^{1/2}\lesi \f{1}{t^{1/2+n/4}}.
	\]
	This implies that, for $y\in B(x_0,r)$,
	\[
	\Big(\int_{B(x_0,2r)\backslash B(x_0,3r/2)} |\nabla_x p_t(x,y)|^2dx\Big)^{1/2}\lesi \f{e^{-\f{r^2}{ct}}}{t^{1/2+n/4}}.
	\]
	Inserting this into \eqref{eq-I2-Dirichlet} we obtain the right hand side of \eqref{eq-lem1D} for $I_1$. This completes our proof.
\end{proof}

\begin{lem}
	\label{lem- vc implies hardy}
	If $f$ is a function supported in a ball $B$ with $4B\cap \partial\Om \ne \emptyset$ and $\|f\|_{L^\vc}\lesi |B|^{-1/p}$, then we have $\|f\|_{H^p_{Mi}(\Om)}\lesi 1$.
\end{lem}
\begin{proof}
	We consider the family of balls $\{B(x,\delta(x)/6): x\in B\}$ which covers the ball $B$. By Vitali's covering lemma we can pick a subfamily of balls denoted by $\{B_j:=B(x_j,\delta(x_j)/2): j\in \mathbb{N}\}$ so that $B\subset \cup_{j\in \mathbb{N}} B(x_j,\delta(x_j)/2)$ and the family $\{\f{1}{3}B_j: j\in \mathbb{N}\}$ are pairwise disjoint. We now write
	\[
	f=\sum_{j}\f{f\chi_{B_j}}{\sum_i \chi_{B_i}}=\sum_{j} \lambda_j A_j
	\]
	where	
	\begin{align*}
	\lambda_j= \Big(\f{|B_j|}{|B|}\Big)^{1/p} \qquad \text{and} \qquad A_j= \Big(\f{|B_j|}{|B|}\Big)^{-1/p}\f{f\chi_{B_j}}{\sum_i \chi_{B_i}}.
\end{align*}
It is clear  that  $A_j$ is an $H^p_{Mi}(\Om)$-atom for every $j$. Indeed note that $\supp A_j\subset B_j$; moreover, we have
	\[
	\|A_j\|_{L^\vc}\le \Big(\f{|B_j|}{|B|}\Big)^{-1/p}\|f\|_{L^\vc}\le \Big(\f{|B_j|}{|B|}\Big)^{-1/p}|B|^{-1/p}=|B_j|^{-1/p}.
	\]
Now since for each $j$ the ball $\f{1}{3}B_j$ is contained in $4B$ then 
	\[
	\sum_j|\lambda_j|^p =\sum_j\f{|B_j|}{|B|}\lesi \sum_j\f{|\f{1}{3}B_j|}{|B|}\lesi 1,
	\]
and this gives $\|f\|_{H^p_{Mi}(\Om)}\lesi 1$.
\end{proof}

We now turn to the proof of Theorem \ref{mainth-Dirichlet}.
\begin{proof}
	[Proof of Theorem \ref{mainth-Dirichlet}:] We shall only give the proof for the case $|\Om|<\vc$. The remaining case $|\Om|=\vc$ can be done in a similar way.
	
	Since Theorem \ref{mainthm1} applies to $L_D$ we may write $H^p_{L_D}$ for any of $H^p_{L_D,\max}$ or $H^p_{L_D,\rad}$. The  inclusion $H^p_{Mi}(\Om)\hookrightarrow H^p_{L_D}(\Om)$ is standard and a similar proof can be found in \cite[Proposition 5.3]{DHMMY}.  Thus we will only demonstrate $H^p_{L_D}(\Om)\hookrightarrow H^p_{Mi}(\Om)$ and to do this we draw upon the atomic characterization in Theorem \ref{mainthm1}. It suffices therefore to prove that for each $(p,M)_{L_D}$-atom $a$ with $M>\f{n}{2}(\f{1}{p}-1)$ we have 
	\begin{align}\label{eq-mi}
	\|a\|_{H^p_{Mi}(\Om)}\lesi 1.
	\end{align}
	If  $a$ is a $(p,M)_{L_D}$-atom of type (a) or (b) from Definition \ref{def: L-atom} associated to a ball $B$, and that ball satisfies $4B\cap \partial\Om\ne \emptyset$, then \eqref{eq-mi} holds by Lemma \ref{lem- vc implies hardy}. It remains to consider the case $4B\subset \Om$. In this case, we have $a=L_D^Mb$. We now write
	\[
	a=L_De^{-r_B^2L_D}\tilde{b}+L_D(I-e^{-r_B^2L_D})\tilde{b}=L_De^{-r_B^2L_D}\tilde{b}+(I-e^{-r_B^2L_D})a=:a_1+a_2
	\]
	where $\tilde{b}=L_D^{M-1}b$. 
	
	We only treat  $a_2$ since $a_1$ can be handled similarly and is easier. To do this let $k_0$ be the positive integer such that $2^{k_0-1}r_B\le \delta(x_B)<2^{k_0}r_B$. Then $k_0\ge 3$ necessarily. We set $S_j(B):=[2^{j+1}B\backslash 2^jB]\cap \Om$ if $j>0$ and $S_0(B):=2B$.We decompose $a_2$ as follows:
	\[
	\begin{aligned}
	a_2&=\sum_{j=k_0-3}^\vc a_2\chi_{S_j(B)}+\sum_{j=0}^{k_0-3}\Big(a_2\chi_{S_j(B)}-\f{\chi_{S_j(B)}}{|S_j(B)|}\int_{S_j(B)}a_2\Big)+\sum_{j=0}^{k_0-3}\f{\chi_{S_j(B)}}{|S_j(B)|}\int_{S_j(B)}a_2\\
	&=:\sum_{j=k_0-3}^\vc 2^{-j}\pi_{1,j}+\sum_{j=0}^{k_0-3}2^{-j}\pi_{2,j}+\sum_{j=0}^{k_0-3}\f{\chi_{S_j(B)}}{|S_j(B)|}\int_{S_j(B)}a_2.
	\end{aligned}
	\]
	For the first summation it is clear that $\supp\,\pi_{1,j}\subset S_j(B)\subset B_j:=2^{j+1}B$ and $4B_j\cap \partial\Om\ne \emptyset$ for all $j\ge k_0-3$. Moreover, we have
	\[
	\pi_{1,j}= 2^j (I-e^{-r_B^2L_D})a\cdot\chi_{S_j(B)}.
	\]
	For $j=0,1,2$ using the $L^\vc$-boundedness of $(I-e^{-r_B^2L_D})$ we have
	\begin{equation}
	\label{eq-j012}
	\|\pi_{1,j}\|_{L^\vc}\lesi \|a\|_{L^\vc}\lesi |B|^{-1/p}\sim |B_j|^{-1/p}.
\end{equation}
	For $j\ge 3$ we use
	\[
	\begin{aligned}
	\pi_{1,j}=2^j\int_0^{r_B^2} sL_De^{-sL_D}a\cdot\chi_{S_j(B)} \f{ds}{s}.
	\end{aligned}
	\]
	and Gaussian bounds on the kernel of $sL_De^{-sL_D}$ (which carry over from (D1)) to obtain
\begin{equation}
\label{eq-j ge 3}
	\begin{aligned}
	\|\pi_{1,j}\|_{L^\vc}&\le 2^j\int_0^{r_B^2} \|sL_De^{-sL_D}a\|_{L^\vc(S_j(B))} \f{ds}{s}\\
	&\lesi 2^j\int_0^{r_B^2} e^{-\f{2^jr_B^2}{cs}}\|a\|_{L^\vc(B)} \f{ds}{s}\\
	&\lesi |2^jB|^{-1/p}=|B_j|^{-1/p}.
	\end{aligned}
	\end{equation}
	From \eqref{eq-j012}, \eqref{eq-j ge 3} and Lemma \ref{lem- vc implies hardy} we have $\|\pi_{1,j}\|_{H^p_{Mi}(\Om)}\lesi 1$, and hence $\sum_{j=k_0-3}^\vc 2^{-j}\pi_{1,j} \in H^p_{Mi}(\Om)$.
	
	For the terms $\pi_{2,j}$, we note that $\int \pi_{2,j} = 0$ and $\supp \pi_{2,j}\subset B_j:=2^{j+1}B$ with $4B_j\subset \Om$. Arguing similarly to the estimates of $\pi_{1,j}$ we also find that $\|\pi_{2,j}\|_{L^\vc}\lesi |B_j|^{-1/p}$. Hence, $\pi_{2,j}$ is an $H^p_{Mi}(\Om)$-atom for each $j$. This implies $\sum_{j=k_0-3}^\vc 2^{-j}\pi_{2,j} \in H^p_{Mi}(\Om)$.
	
	For the last term, we decompose further as follows:
	\[
	\begin{aligned}
		\sum_{j=0}^{k_0-3}\f{\chi_{S_j(B)}}{|S_j(B)|}\int_{S_j(B)}a_2&=\sum_{j=0}^{k_0-3}\Big(\f{\chi_{S_j(B)}}{|S_j(B)|}-\f{\chi_{S_{j-1}(B)}}{|S_{j-1}(B)|}\Big)\int_{2^{k_0-3}\backslash 2^jB}a_2+\f{\chi_{2B}}{|2B|}\int_{2^{k_0-3}B}a_2.
	\end{aligned}
	\]
	Now arguing as above, we can show that for  $j=0,1,\ldots, k_0-3$
	\[
	\Big\|\Big(\f{\chi_{S_j(B)}}{|S_j(B)|}-\f{\chi_{S_{j-1}(B)}}{|S_{j-1}(B)|}\Big)\int_{2^{k_0-3}\backslash 2^jB}a_2\Big\|_{H^p_{Mi}(\Om)}\lesi 2^{-j}.
	\]
	For the remaining term $\f{\chi_{2B}}{|2B|}\int_{2^{k_0-3}B}a_2$  we have
	\[
	\int_{2^{k_0-3}B}a_2=\int_{2^{k_0-3}B}\int_0^{r_B^2} sL_De^{-sL_D}a(x) \f{ds}{s}dx=\int_0^{r_B^2}\int_B \int_{2^{k_0-3}B}q_s(x,y)a(y)dxdy \f{ds}{s}.
	\]
	On the other hand, by Lemma \ref{lem1-Dirichlet} we obtain
	\[
	\Big|\int_{2^{k_0-3}B}q_s(x,y)a(y)dx\Big|\lesi \f{\sqrt{s}}{2^{k_0}r_B}e^{-\f{(2^{k_0}r_B)^2}{cs}}.
	\]
	Hence 
	\[
	\int_{2^{k_0-3}B}a_2\lesi \int_0^{r_B^2}\|a\|_{L^1}\f{\sqrt{s}}{2^{k_0}r_B}e^{-\f{(2^{k_0}r_B)^2}{cs}} \f{ds}{s}\lesi 2^{-k_0}e^{-c2^{2k_0}}|B|^{1-1/p}\lesi |B||2^{k_0-1}B|^{-1/p},
	\]
	which implies
	\[
	\f{\chi_{2B}}{|2B|}\int_{2^{k_0-3}B}a_2\lesi |2^{k_0-2}B|^{-1/p}.
	\]
	As a consequence, $\f{\chi_{2B}}{|2B|}\int_{2^{k_0-3}B}a_2$ is an $H^p_{Mi}(\Om)$-atom associated to the ball $2^{k_0-2}B$. It follows that
	\[
	\Big\|\f{\chi_{2B}}{|2B|}\int_{2^{k_0-3}B}a_2\Big\|_{H^p_{Mi}(\Om)}\lesi 1.
	\]
	This completes our proof.
\end{proof}

\subsection{Schr\"odinger operators with Dirichlet boundary conditions}

Let $\Om$ is a strongly Lipschitz domain of $\Rn$ with $n\ge 3$. This means that $\Om$ is a proper open connected set in $\Rn$ and whose boundary is a finite union of parts of rotated graphs of Lipschitz maps, with at most one of these parts possibly infinite. The class of strongly Lipschitz domains includes special Lipschitz domains, bounded Lipschitz domains and exterior domains. See for example \cite{AR}.

Let $0\le V\in RH_{\tilde{q}}(\Rn)$ with $\tilde{q}>n/2$, i.e.,
\[
\Big(\f{1}{|B|}\int_BV(x)^{\tilde{q}}dx\Big)^{1/\tilde{q}}\le C\f{1}{|B|}\int_BV(x)dx
\]
for all balls $B\subset \Rn$.

We define
\[
W^{1,2}_{V,0}(\Om)=\{u\in W^{1,2}_0(\Om): \int_{\Om}|u(x)|^2V(x)dx<\vc\}.
\]
The Schr\"odinger $L$ on $\Om$ with  Dirichlet Boundary Condition (DBC) can be defined via the following sesquilinear form $\mathcal{Q}$ by setting
\begin{equation}
\label{eq1-defn Schrodinger}
\mathcal{Q}(f,g)=\int_{\Om}\nabla f(x)\overline{\nabla g(x)}dx +\int_{\Om} f(x)\overline{ g(x)}V(x)dx.
\end{equation}
Then $L$ can be written as $L f=-\Delta f +Vf$ where $f\in D(L)$ with
\begin{equation}
\label{eq2-defn Schrodinger}
D(L)=\Big\{f\in W^{1,2}_{V,0}(\Om): \exists g\in L^2(\Om): \mathcal{Q}(f,\phi)=\int_{\Om} g(x)\overline{ \phi(x)}, \ \forall \phi\in  W^{1,2}_{V,0}(\Om)\Big\}.
\end{equation}


For $V\in RH_q, q>n/2$, we define the \emph{critical function} $\rho(x)$ as follows:
\begin{equation}
\label{eq-critical functions}
\rho(x)=\Big\{r\in (0,\vc): \f{1}{r^{n-2}}\int_{B(x_0,r)}V(y)dy\le 1\Big\}.
\end{equation}
Then there exist positive constants $C$ and $k_0$ so that
\begin{equation}\label{criticalfunction}
\ry\leq C\rx\left(1 +\f{d(x,y)}{\rx}\right)^{\f{k_0}{k_0+1}}
\end{equation}
for all $x,y\in X$. See for example \cite{Sh}.

The critical function $\rho$ plays an important role in the rest of this section. Firstly it contributes to better bounds on the heat kernel for $L$ compared to those in (A2) and (A3).

\begin{lem}[\cite{CFYY}]
	\label{lem2-Schrodinger operator}
	Let $L$ be a Schr\"odinger with DBC on the strongly Lipschitz domain with $V\in RH_{\tilde{q}}, \tilde{q}>n/2$. Then we have
	\begin{enumerate}[(i)]
		\item for any $N>0$ there exists $C=C(N)>0$ and $c>0$ so that for all $t>0$ and $x,y\in \Om$,
		\begin{equation}
		\label{eq1- pt Schrodinger}
		0\le p_t(x,y)\le \f{C}{t^{n/2}}\exp\Big(-\f{|x-y|^2}{ct}\Big)\Big[1+\f{\sqrt{t}}{\rho(x)}+\f{\sqrt{t}}{\rho(y)}\Big]^{-N};
		\end{equation}
		\item for any $N>0$ and $0<\delta<\min\{\gamma_0, 2-n/\tilde{q}\}$, there exists $C=C(N,\delta)>0$ and $c>0$ so that for all $t>0$ and $x,y,y'\in \Om$ with $|y-y'|<\sqrt{t}$ and $0<t<{\rm diam}\,\Om$,
		\begin{equation}
		\label{eq1- pt Schrodinger}
		|p_t(x,y)-p_t(x,y')|\le \f{C}{t^{n/2}}\Big(\f{|y-y'|}{\sqrt{t}}\Big)^\delta\exp\Big(-\f{|x-y|^2}{ct}\Big)\Big[1+\f{\sqrt{t}}{\rho(x)}+\f{\sqrt{t}}{\rho(y)}\Big]^{-N};
		\end{equation}
		\item there exist $\alpha>0$ and $C=C(\alpha)$ so that for all $t>0$ and $y\in \Om$,
		\[
		\Big(\int_{\Om}|\nabla_x p_t(x,y)|e^{\f{|x-y|^2}{\alpha t}}dx\Big)^{1/2}\le \f{C}{t^{1/2+n/4}}.
		\]
	\end{enumerate} 
\end{lem}

The function $\rho$ also gives us a useful covering of $\Rn$.
\begin{lem}[\cite{DZ2}]\label{coveringlemm}
	\label{Lem2: rho}
	There exists a family of balls $\{B_\alpha\}_{\alpha\in \mathcal{I}}$ given by $B_\alpha=B(x_\alpha, \rho(x_\alpha))$ satisfies
	\begin{enumerate}[{\rm (i)}]
		\item $\displaystyle \bigcup_{\alpha\in \mathcal{I}} B(x_\alpha, \rho(x_\alpha)) = \Rn$;
		\item For every $\lambda \geq 1$ there exist constants $C$ and $N_1$ such that $\displaystyle \sum_{\alpha\in \mathcal{I}} \chi_{B(x_\alpha, \lambda \rho(x_\alpha))}\leq C\lambda^{N_1}$.
	\end{enumerate}
\end{lem}

Finally the function $\rho$ can be used to define an atomic Hardy space for $L$ which we now present. 
\begin{defn}[Hardy spaces for the DBC Schr\"odinger operator] \label{def-hardyschrodinger}
Let $p\in (\f{n}{n+1},1]$. A bounded, measurable function $a: \Om\to \mathbb{R}$ supported in a ball $B$ is called an $(p,\rho)$-atom if either 
\begin{enumerate}[\upshape(a)]
	\item $a$ is an $H^p_{Mi}$-atom and $r_B< \rho(x_B)/4$; or,
	\item $\|a\|_{L^\vc(\Om)}\le |B|^{-1/p}$, $a\equiv 0$ on $\Om\backslash B$ and  $r_B\ge \rho(x_B)/4$.
\end{enumerate}

We now define the Hardy space $H^{p}_{\rho}(\Om)$ as a set of all $f$ such that
\[
f=\sum_j\lambda_j a_j
\]
where $a_j$ are $(p,\rho)$-atoms and $\lambda_j$ are scalars such that
$
\sum_{j}|\lambda_j|^p<\vc.
$
We also set
\[
\|f\|^p_{H^{p}_{\rho}(\Om)}=\inf\Big\{\sum_{j}|\lambda_j|^p: f=\sum_j\lambda_j a_j \Big\}
\]
where the infimum is taken over all such decompositions.
\end{defn}

\begin{thm}
	\label{mainth-SchrodingerDirichlet}
	Let $p\in (\f{n}{n+\delta},1]$ where $\delta=\min\{\gamma_0,2-n/\tilde{q}\}$. We have
	\[
	H^p_{L,\max}(\Om) \equiv H^p_{L,\rad}(\Om)\equiv H^{p}_{\rho}(\Om), \ \ \f{n}{n+\delta}<p\le 1.
	\]
\end{thm}

\begin{rem} We have the following remarks. 
\begin{enumerate}[(a)]
\item It is important to note that due to Lemma \ref{coveringlemm} we may assume that each $(p,\rho)$-atom satisfying Definition \ref{def-hardyschrodinger} (b)  also satisfies $\rho(x_B)/4<r_B\le \rho(x_B)$.
\item
Just like  the Hardy spaces of Miyachi $H^p_{Mi}(\Om)$, the new Hardy space $H^p_\rho(\Om)$ bears a close relationship with the following Hardy space of restriction related to critical function $\rho$.
For $p\in (\f{n}{n+1},1]$ we say a bounded function $a$ supported in a ball $B\subset \Rn$ is called an $(p,\rho)_{\Rn}$-atom if either 
\begin{enumerate}[(a)]
	\item $\|a\|_{L^\vc(\Om)}\le |B|^{-1/p}$,
	\item $\int a(x)dx=0$ if $r_B<\rho(x_B)/4$.
\end{enumerate}
We define
\[
\|f\|^p_{H^{p}_{\rho}(\Om)}:=\inf\Big\{\sum_{j}|\lambda_j|^p: f=\sum_j\lambda_j a_j \Big\}
\]
where the infimum is taken over all such decompositions $f=\sum_j\lambda_j a_j$ with  $(p,\rho)_{\Rn}$-atoms $a_j$ and numbers $\lambda_j$ satisfying
$\sum_{j}|\lambda_j|^p<\vc$.
Then the Hardy space $H^{p}_{\rho}(\Rn)$ is defined as the completion in the quasi-norm $\|f\|^p_{H^{p}_{\rho}(\Om)}$ of the set
$
\{f\in L^2: f=\sum_j\lambda_j a_j \}
$.

The Hardy space of restriction related to $\rho$ is now defined as
\[
H^p_{\rho, r}(\Om)=\{f: \ \text{there exists $F\in H^p_\rho(\Rn)$ so that $F|_{\Om}=f$}\}
\]
with the norm
\[
\|f\|_{H^p_{\rho,r}(\Om)}=\inf\{\|F\|_{H^p_{\rho}(\Rn)}: F\in H^p_{\rho}(\Rn), F|_{\Om}=f\}.
\]
Then it was proved in \cite{CFYY} that if either $\Om$ is bounded or $\Om^c$ is unbounded, then $H^p_{\rho, r}(\Om)= H^p_{\rho}(\Om)$ for all $\f{n}{n+\delta}<p\le 1$ with $\delta=\min\{\gamma_0,2-n/\tilde{q}\}$. This and Theorem \ref{mainth-SchrodingerDirichlet}  immediately imply the following result:
\begin{cor}
	Let $\Om$ be a strongly Lipschitz domain such that either $\Om$ is bounded or $\Om^c$ is unbounded. Let $p\in (\f{n}{n+\delta},1]$ where $\delta=\min\{\gamma_0,2-n/\tilde{q}\}$. Then we have
	\[
	H^p_{L,\max}(\Om) \equiv H^p_{L,\rad}(\Om)\equiv H^{p}_{\rho,r}(\Om), \ \ \f{n}{n+\delta}<p\le 1.
	\]
\end{cor}
\end{enumerate}
\end{rem}

As in the proof of Theorem \ref{mainth-SchrodingerDirichlet} we require certain kernel estimates first.

\begin{lem}\label{lem1-Schrodinger Dirichlet}
	Let $q_t(x,y)$ be the kernel of $tLe^{-tL}$. Suppose that $x_0\in \Om$ and $0<r<\rho(x_0)/4$ so that $B(x_0,2r)\subset \Om$. Then we have
	\begin{equation}\label{q_t- Schrodinger}
		\Big|\int_{B(x_0,2r)}q_t(x,y)dx\Big| \lesi \Big(\f{\sqrt{t}}{r}\Big)^{2-n/q}
	\end{equation}
	for all $0<t<\rho(x_0)^2$ and all $y\in B(x_0,r)$.
\end{lem}
\begin{proof}
	As in the proof of Lemma \ref{lem1-Dirichlet}, we take $\psi \in C^\vc_c(x_0,2r)$ so that $\psi=1$ on $B(x_0,3r/2)$ and $|\nabla \psi|\lesi 1/r$. Then we have
	\[
	\begin{aligned}
	\int_{B(x_0,2r)}q_t(x,y)dx&=\int_{B(x_0,2r)}q_t(x,y)\psi(x)dx+\int_{B(x_0,2r)}q_t(x,y)[1-\psi(x)]dx
	=:J_1+J_2.
	\end{aligned}
	\]
	We can argue similarly to $I_2$ of Lemma \ref{lem1-Dirichlet} to get $|I_2|\lesssim e^{-r^2/ct}$.
	
	To estimate $I_1$ we use \eqref{eq1-defn Schrodinger} and \eqref{eq2-defn Schrodinger} to deduce
	\begin{equation*}
		\begin{aligned}	
			J_1
			 =t\int_{\Om} \nabla_x p_t(x,y)\cdot \nabla \psi(x)dx+t\int_{\Om} p_t(x,y)V(x)\psi(x)\,dx=:J_{11}+J_{12}.
		\end{aligned}
	\end{equation*}
	Then again arguing as in the proof of estimate $I_1$ from Lemma \ref{lem1-Dirichlet} we can obtain $|J_{11}|\lesssim e^{-r^2/ct}$.
	
	Using \cite[Lemma 5.1]{DZ} we conclude that
	\[
	|J_{12}|\lesi \Big(\f{\sqrt{t}}{\rho(y)}\Big)^{2-n/\tilde{q}}.
	\]
	On the other hand, since $|y-x_0|<r<\rho(x_0)/4$, from \eqref{criticalfunction} we have $\rho(y)\sim \rho(x_0)>4r$. Hence,
	\[
	|J_{12}|\lesi \Big(\f{\sqrt{t}}{r}\Big)^{2-n/\tilde{q}}.
	\]
	Collecting all estimates $I_2, J_{11}$ and $J_{12}$ we get the desired estimate \eqref{q_t- Schrodinger}.
\end{proof}

We are now ready to give the proof of Theorem \ref{mainth-SchrodingerDirichlet}.

\begin{proof}[Proof of Theorem \ref{mainth-SchrodingerDirichlet}:] 
	We shall only give the proof for the case $|\Om|<\vc$ since the remaining case $|\Om|=\vc$ can be done similarly.

	We will first show that $H^p_\rho(\Om)\subset H^p_{L,\rad}(\Om)$. Indeed, let $a$ be a $(p,\rho)$-atom associated to a ball $B$. Now if $r_B<\rho(x_B)/4$, then in this case $a$ is a $H^p_{Mi}(\Om)$-atom and a standard argument shows that
	\[
	\left\|a^+_L\right\|_{L^p(\Om)}\lesi 1.
	\]
On the other hand if $\rho(x_B)/4<r_B\le \rho(x_B)$, then we split
$$
\|a^+_L\|^p_{L^p}\leq \|a^+_L\|^p_{L^p(4B)}+\|a^+_L\|^p_{L^p(\Omega\backslash 4B)}: =I_1+I_2.
$$
H\"older's inequality and the estimate $\|a^+_L\|_{L^\vc}\lesi \|a\|_{L^\vc}$ allow us to  readily conclude that $I_1\lesi 1$.

We turn to the second term $I_2$. Note firstly that $r_B\sim \rho(x_0)\sim \rho(y)$, and secondly that $|x-y|\sim |x-x_0|$ holds whenever $y,x_0\in B$ and $x\in (4B)^c$. These facts in tandem with  \eqref{eq1- pt Schrodinger}
allow us to obtain, for $N>n(1-p)/p$,
$$
\begin{aligned}
I_2&\lesi \int_{(4B)^c}\sup_{t>0}\Big[\int_B \f{1}{t^{n/2}}\exp\Big(-\f{|x-y|^2}{ct}\Big)\Big(\f{\rho(y)}{\sqrt{t}}\Big)^N|a(y)|dy\Big]^pdx\\
%
&\lesi \int_{(4B)^c}\sup_{t>0}\Big[\int_B \f{1}{t^{n/2}}\exp\Big(-\f{|x-x_0|^2}{ct}\Big)\Big(\f{r_B}{\sqrt{t}}\Big)^N|a(y)|dy\Big]^pd\mu(x)\\
&\lesi \int_{(4B)^c}\Big[\int_B \f{1}{|x-x_0|^n}\Big(\f{r_B}{|x-x_0|}\Big)^N|a(y)|dy\Big]^pdx\\
&\lesi 1.
\end{aligned}
$$
This completes the direction $H^p_\rho(\Om)\subset H^p_{L,\rad}(\Om)$. The reverse direction can be done in a similar way to that of Theorem \ref{mainth-Dirichlet} and will be omitted.	 
\end{proof}

\subsection{Fourier--Bessel operators on $((0,1),dx)$}

For $\nu>-1$, we consider the following differential operator

\begin{align}\label{eq-FBdx}
L=-\f{d^2}{dx^2} +\f{\nu^2-1/4}{x^2}.
\end{align}
Let $\{\lambda_{k,\nu}\}_{k\ge 1}$ denote the sequence of successive positive zeros of the Bessel function $J_\nu$ and consider 
$$
\psi_k^\nu(x)=d_{k,\nu}\lambda^{1/2}_{k,\nu}J_\nu(\lambda_{k,\nu}x)x^{1/2}
$$
where $x\in (0,1)$ and $d_{k,\nu}=\sqrt{2}|\lambda_{k,\nu}J_{\nu+1}(\lambda_{k,\nu})|^{-1}$.

Then the  system $\{\psi_k^\nu\}_k$ forms an orthornomal basis for $((0,1),dx)$. It is well-known that
\[
 L\psi_k^\nu(x)=\lambda_{k,\nu}^2\psi_k^\nu(x).
\]
The operator $L$ has a non-negative self-adjoint extension which is still denoted by $L$ with domain
\[
D(L)=\{f\in L^2((0,1),dx): \sum_{k=1}^\infty \lambda_{k,\nu}^4 |\langle f, \psi_k^\nu\rangle|^2<\infty\}.
\]
This operator is called the Bessel operator on  $((0,1),dx)$.

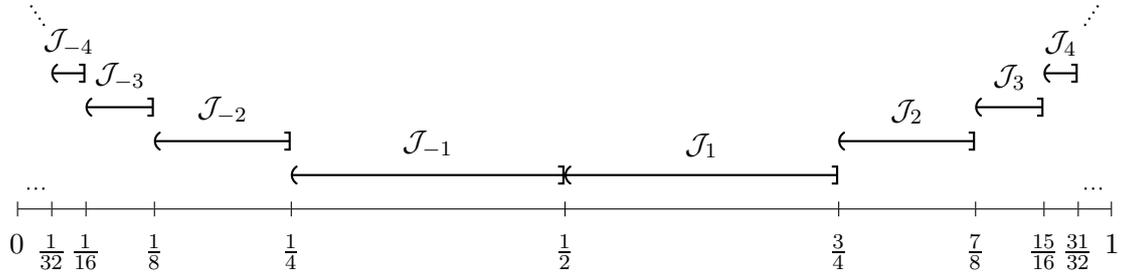
\begin{figure}[h]
\begin{center}
\begin{tikzpicture}[scale=0.9]
\path (-4,-1) coordinate (L1);
\path (-6,-1) coordinate (L2);
\path (-7,-1) coordinate (L3);
\path (-7.5,-1) coordinate (L4);
\path (-8,-1) coordinate (L5);
\path (0,-1) coordinate (00);
\path (4,-1) coordinate (R1);
\path (6,-1) coordinate (R2);
\path (7,-1) coordinate (R3);
\path (7.5,-1) coordinate (R4);
\path (8,-1) coordinate (R5);

\draw [|-] (L5)--(L4);
\draw [|-] (L4)--(L3);
\draw [|-] (L3)--(L2);
\draw [|-] (L2)--(L1);
\draw [|-] (L1)--(00);
\draw [|-] (00)--(R1);
\draw [|-] (R1)--(R2);
\draw [|-] (R2)--(R3);
\draw [|-] (R3)--(R4);
\draw [|-|] (R4)--(R5);

\draw[[-),thick] (4,-.5) -- (0,-.5);
\draw[[-),thick] (6,0) -- (4,0);
\draw[[-),thick] (7,0.5) -- (6,0.5);
\draw[[-),thick] (7.5,1) -- (7,1);
\draw[[-),thick] (0,-.5) -- (-4,-.5);
\draw[[-),thick] (-4,0) -- (-6,0);
\draw[[-),thick] (-6,0.5) -- (-7,0.5);
\draw[[-),thick] (-7,1) -- (-7.5,1);

\node[above=1mm] at (2,-.5) {$\mathcal{J}_1$};
\node[above=1mm] at (5,0) {$\mathcal{J}_2$};
\node[above=1mm] at (6.5,.5) {$\mathcal{J}_3$};
\node[above=1mm] at (7.25,1) {$\mathcal{J}_4$};
\node[above=1mm] at (-2,-.5) {$\mathcal{J}_{-1}$};
\node[above=1mm] at (-5,0) {$\mathcal{J}_{-2}$};
\node[above=1mm] at (-6.5,.5) {$\mathcal{J}_{-3}$};
\node[above=1mm] at (-7.25,1) {$\mathcal{J}_{-4}$};
\node[below=2mm] at (00) {$\f{1}{2}$};
\node[below=2mm] at (R1) {$\f{3}{4}$};
\node[below=2mm] at (R2) {$\f{7}{8}$};
\node[below=2mm] at (R3) {$\f{15}{16}$};
\node[below=2mm] at (R4) {$\f{31}{32}$};
\node[below=2mm] at (R5) {$1$};
\node[below=2mm] at (L1) {$\f{1}{4}$};
\node[below=2mm] at (L2) {$\f{1}{8}$};
\node[below=2mm] at (L3) {$\f{1}{16}$};
\node[below=2mm] at (L4) {$\f{1}{32}$};
\node[below=2mm] at (L5) {$0$};

\draw[dotted, thick] (7.6,1.7) -- (7.8,2);
\draw[dotted, thick] (-7.6,1.7) -- (-7.8,2);
\draw[dotted, thick] (-7.6,-.7)--(-7.9,-.7);
\draw[dotted, thick] (7.6,-.7)--(7.9,-.7);

 \end{tikzpicture}
\end{center}
\caption{Intervals $\mathcal{J}_j$} \label{fig: dx}
\end{figure}

In order to consider the maximal function characterization for the Hardy spaces associated to $L$, as in \cite{BDK} we consider the intervals:
\begin{equation}\label{eq-Jj}
\mathcal{J}_j=\begin{cases}
(1-2^{-j},1-2^{-j-1}], \ \ \ &j\ge 1\\
(2^{j-1},2^j], & j\le -1.
\end{cases}
\end{equation}
which are depicted in Figure \ref{fig: dx}. 

It is obvious that the family $\{\mathcal{J}\}_{j\in \mathbb{N}}$ is pairwise disjoint and $(0,1)=\bigcup_{j\in \mathbb{N}}\mathcal{J}_j$. For each $j\in \mathbb{Z}^*$ we also denote $\mathcal{J}_j^*=\tfrac{1}{10}\mathcal{J}_j$. We now consider the following atoms.
\begin{defn}\label{def: atoms 1}
	Let $\f{1}{2}<p\le 1$. A function $a$ is a type (a) atom associated with an interval $I\subseteq (0,1)$ if 
	\begin{enumerate}[\upshape(i)]
		\item $\supp a\subset I$
		\item $\Vert a\Vert_{L^\vc}\le |I|^{-1/p}$
		\item $\int a(x)dx=0$
	\end{enumerate}
	A function $a$ is a type (b) atom if
	$$a(x)=\f{\chi_{\mathcal{J}_j}(x)}{|\mathcal{J}_j|^{1/p}}$$ for some $j\in \mathbb{N}$.
	We say a function $a$  is  an $H^p((0,1),dx)$-atom  if it is either a type (a) or type (b) atom.
\end{defn}

\begin{defn}[Atomic Hardy spaces on $((0,1),dx)$]
	Given $p\in (\f{1}{2},1]$, $q\in [1,\vc]\cap (p,\vc]$, we  say that $f=\sum
	\lambda_ja_j$ is an atomic $H^p((0,1),dx)$-representation if
	$\{\lambda_j\}_{j=0}^\infty\in l^p$, each $a_j$ is an $H^p((0,1),dx)$-atom,
	and the sum converges in $L^2(X)$. The space $H^{p}_{at}((0,1),dx)$ is then defined as the completion of
	\[
	\left\{f\in L^2:f \ \text{has an atomic
		$H^p((0,1),dx)$-representation}\right\},
	\]
	with the norm given by
	$$
	\|f\|_{H^{p}_{at}((0,1),dx)}=\inf\left\{\left(\sum|\lambda_j|^p\right)^{1/p}:
	f=\sum \lambda_ja_j \ \text{is an atomic $H^p((0,1),dx)$-representation}\right\}.
	$$
\end{defn}
As in  Definition \ref{defn-maximal Hardy spaces} we denote by $H^p_{L, \rad}((0,1),dx)$ and $H^p_{L, \max}((0,1),dx)$ respectively the maximal Hardy spaces defined via radial and non-tangential maximal functions associated to $L$.

Our main result in this section is the following:

\begin{thm}\label{Th1-Bessel 1} 
	Let $\nu>-1/2$. Let $p\in (\f{1}{1+\delta},1]$ where $\delta=\min\{1,\nu+1/2\}$. Then we have
	\[
	H^{p}_{at}((0,1),dx)\equiv H^p_{L, \rad}((0,1),dx)\equiv H^p_{L, \max}((0,1),dx)
	\]
	with equivalent norms.
\end{thm}
\begin{rem}
Note that it was proved in \cite{DPRS} that  
\[
H^{1}_{at}((0,1),dx)\equiv H^1_{\sqrt{L}, \rad}((0,1),dx)\equiv H^1_{\sqrt{L}, \max}((0,1),dx).
\]
Hence the results of Theorem \ref{Th1-Bessel 1} can be viewed as an extension of those in  \cite{DPRS} to the range $p$ below 1. 
\end{rem}

In order to give the proof of Theorem \ref{Th1-Bessel 1} we need the following technical material. Firstly we collect together some estimates on the kernels and their derivatives. 
\begin{lem}[\cite{MSZ} in Theorem 1.1]\label{lem-heat kernel-Bessel 1}
	For $\nu>-1$ we have 
	$$p_t(x,y) \approx \f{(xy)^{\nu+1/2}(1+t)^{\nu+2}}{(t+xy)^{\nu+1/2}} \Big(1\land \f{(1-x)(1-y)}{t}\Big)\f{1}{\sqrt{t}}e^{-\f{|x-y|^2}{4t}-\lambda^2_{1,\nu}t} $$
	for $x,y\in X$ and $t>0$.
\end{lem}
When $\nu>-1/2$, a simple calculation applied to the bounds in Lemma \ref{lem-heat kernel-Bessel 1} gives
\begin{equation}\label{eq-GU pxy-Bessel 1}
|p_t(x,y)|\lesi \f{(xy)^{\nu+1/2}}{(t+xy)^{\nu+1/2}}\f{1}{\sqrt{t}}\Big(1\land \f{(1-x)(1-y)}{t}\Big)e^{-\f{|x-y|^2}{ct}}
\end{equation}
for all $x,y\in X$ and $t>0$.
Note that when $-1<\nu<-1/2$ the Gaussian upper bounds for the kernel $p_t(x,y)$ may  fail. Hence and for this reason we restrict ourself to the case $\nu>-1/2$.

\begin{lem}\label{lem: heatderiv 1 -Bessel 1}
	For $\nu>-1$ we have
	\begin{align}\label{heatderiv1-Bessel 1}
	|\partial_x p_t(x,y)| \lesssim \f{1}{t}e^{-(x-y)^2/4t}+\f{1}{x}p_t(x,y).
	\end{align} 
	for all $x,y\in X$ and $t>0$.
	
	As a consequence, if $\nu>-1/2$, then we have
	\begin{align}\label{heatderiv2-Bessel 1}
	|\partial_x p_t(x,y)| \lesssim \f{1}{xt}e^{-(x-y)^2/4t}
	\end{align} 
		for all $x,y\in X$ and $t>0$.
\end{lem}

\begin{proof}
	We have
	\[
	|\partial_x p_t(x,y)|\lesi \big|\big(\partial_x - \tfrac{\nu+1/2}{x}\big) p_t(x,y) \big|+\f{|\nu+1/2|}{x}p_t(x,y).
	\]
	Using the argument in \cite[Lemma 2.4]{DPRS} we find that
	\[
	\big|\big(\partial_x - \tfrac{\nu+1/2}{x}\big) p_t(x,y) \big| \lesssim \f{1}{t}e^{-(x-y)^2/4t}.
	\]
	This yields \eqref{heatderiv1-Bessel 1} as desired.
	
	The estimate \eqref{heatderiv2-Bessel 1} follows from \eqref{heatderiv1-Bessel 1} and \eqref{eq-GU pxy-Bessel 1}.
\end{proof}

Let us define the notion of ``intervals'' in $(0,1)$. For $x\in (0,1)$ and $r>0$ we denote  by
\begin{align}\label{eq-interval}
I_r(x)=(x-r,x+r)\cap (0,1)
\end{align}
the interval centred at $x$ of radius $r$. Henceforth and unless otherwise specified, by an interval $I$ in $(0,1)$ we shall mean $I=I_{r_I}(x_I)$  for some fixed centre $x_I$ and radius $r_I$. 

We define the \emph{critical function} for $L$ by 
\begin{equation}\label{eq-critical function-B1}
\rho(x) :=\f{1}{3}\min\{x,1-x\}, \qquad x\in (0,1)
\end{equation}

 For $x\in (0,1)$ and $\rho$ defined as in \eqref{eq-critical function-B1}, we denote $I_\rho(x)=I_{\rho(x)}(x)$. Such an interval is called a \emph{critical interval}. 

We have the following result whose easy proof we omit. 

\begin{lem}
	\label{lem-critical function}
	If $I$ is an interval with $r_I\le \rho(x_I)$ then we have, for all $x\in I$
	\begin{enumerate}[{\rm (i)}]
	\item $x\sim x_I$;
	\item $\rho(x)\sim\rho(x_I)$;
	\item $r_I\lesi x$.
\end{enumerate}
\end{lem}
Denote by $Q_t$ the operator
$$ Q_t:=tL e^{-tL}$$
and $q_t(x,y)$ the kernel of $Q_t$. It is well-known that the Gaussian upper bound can be transfered to the kernel $q_t(x,y)$, i.e.,
\begin{equation}\label{eq-GU qxy-Bessel 1}
|q_t(x,y)|\lesi \f{1}{\sqrt{t}}e^{-\f{|x-y|^2}{ct}}
\end{equation}
for all $x,y\in X$ and $t>0$.
We apply \eqref{eq-GU qxy-Bessel 1} to obtain the following.
\begin{lem}\label{lem: Qest dx}
	For any interval $I$ with $r_I\le \rho(x_I)$ we have
	\begin{align*}
	\Big|\int_{I}q_t(x,y)dx\Big| \lesi  \f{t}{r_I^2}
	\end{align*}
	for any $\forall y\in\tfrac{1}{2}I$ and $t>0$. 
\end{lem}

\begin{proof}
	Define the cutoff function $\varphi\in C^\infty_c(X)$ supported in $I$ with $0\le \varphi\le 1$, equal to 1 on $\f{3}{4}I$ and whose derivative satisfies $|\varphi'(x)|\lesssim 1/r_I$.
	Then 
	\begin{align*}
	\Big|\int_{I}q_t(x,y)dx\Big|&\lesi  \Big|\int_{X}t\partial_tp_t(x,y)\varphi(x)dx\Big|+\Big|\int_{I\backslash \tfrac{3}{4}I }q_t(x,y)[1-\varphi(x)]dx\Big|
	=:I_1+I_2.
	\end{align*}
	For the term $I_2$, using \eqref{eq-GU qxy-Bessel 1} and that $|x-y|\sim r_I$ we have
	\[
	\begin{aligned}
	I_2&\lesi \int_{I\backslash \tfrac{3}{4}I }\f{1}{\sqrt{t}}e^{-\f{|x-y|^2}{4t}}dx
	\lesi e^{-\f{r_I^2}{ct}}\int_{I\backslash \tfrac{3}{4}I }\f{1}{\sqrt{t}}e^{-\f{|x-y|^2}{8t}}dx 
	\lesi e^{-\f{r_I^2}{ct}}.
	\end{aligned}
	\]
	For the first term,  since $\partial_t p_t(\cdot,y)=-Lp_t(\cdot,y)$, then 
	\[
	\begin{aligned}
	I_1
	&\lesi\Big|\int_{X}t\partial^2_{xx}p_t(x,y)\varphi(x)dx\Big|+|\nu^2-1/4|\Big|\int_{\f{3}{4}I}tp_t(x,y)\varphi(x)\f{dx}{x^2}\Big|
	=:I_{11}+I_{12}.
	\end{aligned}
	\]
	Now \eqref{eq-GU pxy-Bessel 1} and (iii) in Lemma \ref{lem-critical function} implies that
$I_{12}\lesi t/r_I^2$.

	Integration by parts gives
	\[
	I_{11}=\Big|\int_{X}t\partial_{x}p_t(x,y)\partial_{x}\varphi(x)dx\Big|,
	\]
	and along with estimate \eqref{heatderiv2-Bessel 1}, part (iii) of Lemma \ref{lem-critical function} and the fact that $|\varphi'(x)|\lesi r_I^{-1}$ yields
	\[
	\begin{aligned}
	I_{11}\lesi \f{t}{r_I}\Big|\int_{\f{3}{4}I}\f{1}{x\sqrt{t}}e^{-\f{|x-y|^2}{4t}}dx\Big|\sim \f{t}{r_I}\Big|\int_{\f{3}{4}I}\f{1}{r_I\sqrt{t}}e^{-\f{|x-y|^2}{4t}}dx\Big|
	\lesi \f{t}{r_I^2},
	\end{aligned}
	\]
	completing our proof.
\end{proof}

We now turn to the action of the radial maximal operator on atoms. Note that the intervals $\mathcal{J}_j^*$ has been defined in the comments after \eqref{eq-Jj}.

\begin{lem}\label{lem1- atom type ab}
Let $\nu>-1/2$  and $\f{1}{1+\delta}<p\le 1$ with $\delta=\min\{1,\nu+1/2\}$.  Suppose that $a$ is either \begin{enumerate}[\upshape(i)]
	\item  a type (b) $H^p\big((0,1),dx\big)$-atom, or 
	\item  a type (a) $H^p((0,1),dx)$-atom supported in $\mathcal{J}_j^*$ for some $j\in \mathbb{Z}^*$.
	\end{enumerate}
Then there exists $C>0$ independent of $a$ so that
	\[
	\Big\|\sup_{t>0}|e^{-tL}a|\Big\|_{L^p((0,1),dx)}\le C.
	\]
\end{lem}

\begin{proof}[Proof of Lemma \ref{lem1- atom type ab}]
Proof of part (i). 
	Since $a$ is an 	$H^p((0,1),dx)$-atom of type (b), then $a=\f{\chi_{I}}{|I|^{1/p}}$ where $I\equiv \mathcal{J}_j$ some $j\in \mathbb{Z}\backslash\{0\}$.
	\[
	\Big\|\sup_{t>0}|e^{-tL}a|\Big\|_{L^p((0,1),dx)}\lesi \Big\|\sup_{t>0}|e^{-tL}a|\Big\|_{L^p(2I)}+\Big\|\sup_{t>0}|e^{-tL}a|\Big\|_{L^p((2I)^c)}=:E_1+E_2.
	\]
	It is easy to see that
	\[
	E_1\lesi |2I|^{1/p}\Big\|\sup_{t>0}|e^{-tL}a|\Big\|_{L^\vc(2I)}\lesi |2I|^{1/p}\|a\|_{L^\vc(2I)}=1.
	\]
We handle $E_2$ by studying the pointwise bounds on $\sup_{t>0}|e^{-tL}a(x)|$. Firstly by the heat kernel bounds \eqref{eq-GU pxy-Bessel 1}, and that $|x-y|\sim |x-x_I|$ whenever $x\in (2I)^c$, we have
\[
\begin{aligned}
\sup_{t>0}|e^{-tL}a(x)|&\lesi \sup_{0<t<r^2_{I}}|I|^{-1/p}\int_{I}\f{1}{\sqrt{t}}e^{-\f{|x-y|^2}{ct}}dy\\
& \qquad +\sup_{t\ge r^2_{I}}|I|^{-1/p}\int_{I}\f{(xy)^{\nu+1/2}}{(t+xy)^{\nu+1/2}}\f{1}{\sqrt{t}}\Big(1\land \f{(1-x)(1-y)}{t}\Big)e^{-\f{|x-x_{I}|^2}{ct}}dy\\
&=: E_{21}(x)+E_{22}(x).
\end{aligned}
\]	
It is straightforward that
\[
E_{21}(x)
\lesi |I|^{1-1/p}\f{1}{|x-x_{I}|}\f{r^2_I}{|x-x_I|^2}
\]	
which implies $\|E_{21}\|_{L^p((2I)^c)}\lesi 1$ provided $p\in (1/2,1)$.

We divide the calculation for $E_{22}$ into two cases.

\textbf{Case 1: $I\equiv \mathcal{J}_j, j>0$.} In this case we have $(1-x)\lesi |x-x_I|$ and $(1-y)\sim r_I$. Hence
\[
\begin{aligned}
E_{22}(x)&\lesi  \sup_{t\ge r^2_{I}}|I|^{-1/p}\int_{I}\f{1}{\sqrt{t}}\f{|x-x_I|r_I}{t}e^{-\f{|x-x_{I}|^2}{ct}}dy
\lesi  |I|^{1-1/p}\f{1}{|x-x_I|}\f{r_I}{|x-x_I|},
\end{aligned}
\]
which implies $\|E_{22}\|_{L^p((2I)^c)}\lesi 1$ whenever $p\in (1/2,1)$.

\textbf{Case 2: $I\equiv \mathcal{J}_j, j<0$.} In this case we have $x\lesi |x-x_I|$ and $y\sim r_I$. Hence,
\[
\begin{aligned}
E_{22}(x)&\lesi \sup_{t\ge r^2_{I}}|I|^{-1/p}\int_{I}\Big(\f{xy}{t}\Big)^{\nu+1/2}\f{1}{\sqrt{t}}e^{-\f{|x-x_{I}|^2}{ct}}dy\\
&\lesi  \sup_{t\ge r^2_{I}}|I|^{-1/p}\int_{I}\Big(\f{r_I|x-x_I|}{t}\Big)^{\nu+1/2}\f{1}{\sqrt{t}}e^{-\f{|x-x_{I}|^2}{ct}}dydy\\
&\lesi |I|^{1-1/p}\f{1}{|x-x_I|}\Big(\f{r_I}{|x-x_I|}\Big)^{\nu+1/2}
\end{aligned}
\]
which yields  $\|E_{22}\|_{L^p((2I)^c)}\lesi 1$, provided that $p\in (\f{1}{1+\delta},1]$ with $\delta=\min\{1, \nu+1/2\}$. 

Collecting together the estimates for $E_{21}$ and $E_{22}$ we obtain $E_2\lesssim 1$, completing the proof of part (i). 

%

	We now prove part (ii). 
	Suppose that $a$ is an $H^p((0,1),dx)$-atom of type (a) associated to some interval $I\subset \mathcal{J}_j^*$ . 
	
	We write
	\[
	\Big\|\sup_{t>0}|e^{-tL}a|\Big\|_{L^p((0,1),dx)}\lesi \Big\|\sup_{t>0}|e^{-tL}a|\Big\|_{L^p(2I)}+\Big\|\sup_{t>0}|e^{-tL}a|\Big\|_{L^p((2I)^c)}=:F_1+F_2.
	\]
 By  arguing similarly to $E_1$ in the proof of part (i)  we have $F_1\lesi 1$.
	
	To handle $F_2$ we use the cancellation property of $a$ to write
	\begin{align*}
	e^{-tL}a(x)=\int_I [p_t(x,y)-p_t(x,x_I)]a(y)\,dy.
	\end{align*}	
	Then for $x\in (2I)^c$ we may apply Lemma \ref{lem: heatderiv 1 -Bessel 1}, the bounds \eqref{eq-GU pxy-Bessel 1}, and the fact that $|x-y|\sim |x-x_I|$  whenever $y\in I$ to obtain
	\[
	\begin{aligned}
	\sup_{t>0}|e^{-tL}a(x)|
	&\lesi \sup_{t>0}\int_I\f{|y-x_I|}{\sqrt{t}}\f{1}{\sqrt{t}}e^{-\f{|x-x_I|^2}{ct}}|a(y)|dy\\
	& \qquad +\sup_{t>0}\int_I\Big(\f{xy}{t+xy}\Big)^{\nu+1/2}\f{|y-x_I|}{y}\f{1}{\sqrt{t}}e^{-\f{|x-x_I|^2}{ct}}|a(y)|dy\\
	&=: F_{21}(x)+F_{22}(x).
	\end{aligned}
	\]
	Since the variable $y$ belongs to $I$ it is  then clear that
	\[
	\begin{aligned}
	F_{21}(x)&\lesi \sup_{t>0}\|a\|_{L^\vc}\int_I\f{r_I}{\sqrt{t}}\f{1}{\sqrt{t}}e^{-\f{|x-x_I|^2}{ct}}dy
	\lesi |I|^{1-1/p}	\f{r_I}{|x-x_I|}\f{1}{|x-x_I|}
	\end{aligned}
		\]
	which implies  $\|F_{21}\|_{L^p((2I)^c)}\lesi 1$.
	
	For the expression $F_{22}$ we further subdivide
	$$ F_{22} = F_{22}\,\chi_{3\JJ_j \backslash 2I} + F_{22}\,\chi_{(3\JJ_j)^c}.$$
		
	Now whenever $x\in 3\JJ_j\backslash 2I$ we have the inequality $y^{-1}\lesi |x-x_I|^{-1}$. Thus the first term can be controlled by
		\[
	\begin{aligned}
		F_{22}(x)\,\chi_{3\JJ_j \backslash 2I}(x)&\lesi \sup_{t>0}\int_I\f{|y-x_I|}{|x-x_I|}\f{1}{\sqrt{t}}e^{-\f{|x-x_I|^2}{ct}}|a(y)|\,dy
		\lesi |I|^{1-1/p}\f{r_I}{|x-x_I|}\f{1}{|x-x_I|},
	\end{aligned}
	\]
	which yields $\|F_{22}\|_{L^p(3\JJ_j\backslash 2I)}\lesi 1$.
	
	For the second term we consider two cases.
	
	\textbf{Case 1: $ j>0$}. In this situation $y \sim 1$, implying $y^{-1}\lesi |x-x_I|^{-1}$  and therefore,
	\[
	\begin{aligned}
	F_{22}(x)\,\chi_{(3\JJ_j)^c}(x)&\lesi \sup_{t>0}\int_I\f{r_I}{|x-x_I|}\f{1}{\sqrt{t}}e^{-\f{|x-x_I|^2}{ct}}|a(y)|dy
	\lesi |I|^{1-1/p}\f{r_I}{|x-x_I|}\f{1}{|x-x_I|}	.
	\end{aligned}
	\]
	
	\textbf{Case 2:  $ j<0$}. In this case $(3\JJ_j)^c=(6r_{\JJ_j},1)$ and hence $x\sim x-x_I$. Then we have
	\[
	\begin{aligned}
	F_{22}(x)\,\chi_{(3\JJ_j)^c}(x)&\lesi \sup_{t>0}\int_I\Big(\f{xy}{t}\Big)^{\delta}\f{r_I}{y}\f{1}{\sqrt{t}}e^{-\f{|x-x_I|^2}{ct}}|a(y)|dy\\
	&\lesi \sup_{t>0}\int_I\Big(\f{y\,|x-x_I|}{t}\Big)^{\delta}\Big(\f{r_I}{y}\Big)^{\delta}\f{1}{|x-x_I|}e^{-\f{|x-x_I|^2}{2ct}}|a(y)|dy\\
	&\lesi |I|^{1-1/p}\Big(\f{r_I}{|x-x_I|}\Big)^{\delta}\f{1}{|x-x_I|}	
	\end{aligned}
	\] 
	where $\delta=\min\{1,\nu+1/2\}$.
	
	Taking into account the bounds in both cases we  conclude $\|F_{22}\|_{L^p((3\JJ_j)^c)}\lesi 1$. 
	
	On combining our estimates for $F_{21}$ and $F_{22}$ we then have $F_2\lesssim 1$, completing our proof of the Lemma. 
%
%
%
%
%
\end{proof}

\medskip

We are now ready to prove the main theorem of this section. 

\begin{proof}[Proof of Theorem \ref{Th1-Bessel 1}:] We split the proof into two steps.

	\textbf{Step 1: $H^p_{L,\rad}((0,1),dx)\subset H^p_{at}((0,1),dx)$.} 
	
	Suppose that $a$ is a $(p,M)_L$-atom as in Definition \ref{def: L-atom} associated to an interval $I$. We consider two cases: $4I\cap (0,1)^c\ne \emptyset$ and $4I\subset (0,1)$. 
	
	\noindent\textbf{Case 1: $4I\cap (0,1)^c\ne \emptyset$.} In this situation, it easy to see that if $x_I\in \mathcal{J}_j$ for some $j\in \mathbb{Z}$, then  $|I|\sim |\mathcal{J}_j|$. Hence, using the decomposition
	\[
	a=\Big[a-\f{\chi_{\mathcal{I}_j}}{|\mathcal{I}_j|}\int a\Big] + \f{\chi_{\mathcal{I}_j}}{|\mathcal{I}_j|}\int a=:\tilde{a}_1+\tilde{a}_2.
	\]
	We see that $a_1$ is an $H^p((0,1),dx)$ of type (a), while $a_2$ is an $H^p((0,1),dx)$ atom of type (b). Thus $a\in H^p_{at}((0,1),dx)$.

	\medskip
	
	\noindent\textbf{Case 2: $4I\subset (0,1)$.} In this case, $a$ can be expressed in the form $a=Lb$. We now write
	\[
	a=Le^{-r_I^2L}b+L(I-e^{-r_I^2L})b=Le^{-r_I^2L}b+(I-e^{-r_I^2L})a=a_1+a_2
	\]
	where $b$ is supported in $B$ and satisfies
	$
	\|b\|_{L^\vc}\le r_I^2 |I|^{-1/p}
	$.
	
	We take care $a_2$ only, since $a_1$ can be similarly treated. We  choose $k_0\in \mathbb{N}$ so that $2^{k_0-1}r_I\le \f{4}{3}\min\{x_I,1-x_I\}=4\rho(x_I) <2^{k_0}r_I$. Hence $k_0\ge 3$. We set $S_j(I)=[2^{j+1}I\backslash 2^jI]\cap (0,1)$ if $j>0$ and $S_0(I)=2I$. Then as in the proof of Theorem \ref{mainth-Dirichlet} we decompose $a_2$ as follows:
	\[
	\begin{aligned}
	a_2=\ &\sum_{j=k_0-3}^\vc a_2\chi_{S_j(I)}+\sum_{j=0}^{k_0-3}\Big(a_2\chi_{S_j(I)}-\f{\chi_{S_j(I)}}{|S_j(I)|}\int_{S_j(I)}a_2\Big)\\
	&+\sum_{j=0}^{k_0-3}\Big(\f{\chi_{S_j(I)}}{|S_j(I)|}-\f{\chi_{S_{j-1}(I)}}{|S_{j-1}(I)|}\Big)\int_{2^{k_0-3}\backslash 2^jI}a_2+\f{\chi_{2I}}{|2I|}\int_{2^{k_0-3}I}a_2\\
	=:\ & A_1 + A_2 + A_3 + A_4.
	\end{aligned}
	\]
	By arguing in a similar way to the proof of Theorem \ref{mainth-Dirichlet} we can show that firstly $A_1$ can be expressed as an atomic representation of type (b) atoms of Definition \ref{def: atoms 1}; and  secondly that $A_2$ and $A_3$ can be expressed as an atomic representation of type (a) atoms.

	It remains then to take care of $A_4$. Firstly note that $\supp A_4\subset 2^{k_0}I=:\widehat{I}$. Next recall that $q_s(x,y)$ is the kernel of $sLe^{-sL}$. Then applying Lemma \ref{lem: Qest dx} we have
	\[
	\begin{aligned}
	\int_{2^{k_0-3}I}a_2&=\int_{2^{k_0-3}I}\int_0^{r_I^2} sLe^{-sL}a(x) \f{ds}{s}dx
	\lesi \int_0^{r_I^2}\f{s}{(2^{k_0-3}r_I)^2}\int_I |a(y)| dy \f{ds}{s}
	\lesi 2^{-2k_0}|I|^{1-1/p},
	\end{aligned}
	\]
and since $p>1/2$ then 
	\[
	\|A_4\|_{L^\vc}\lesi \f{2^{-2k_0}}{|I|^{1/p}}\lesi \f{1}{|2^{k_0}I|^{1/p}} =|\widehat{I}|^{-1/p}.
	\]
	Now since $4\widehat{I}\cap (0,1)^c\ne \emptyset$ we may evoke Case 1 to obtain $\|A_4\|_{H^p_{at}((0,1),dx)}\lesi 1$.

	\textbf{Step 2: $H^p_{at}((0,1),dx)\subset H^p_{L,\rad}((0,1),dx)$.} It suffices to prove that there exists $C>0$ so that 
	\begin{equation}\label{eq-proof step 2}
	\Big\|\sup_{t>0}|e^{-tL}a|\Big\|_{L^p((0,1),dx)}\le C
	\end{equation}
	for all $H^p((0,1),dx)$-atoms $a$. Now if $a$ is type (b) atom then the inequality \eqref{eq-proof step 2} follows from part (i) of Lemma \ref{lem1- atom type ab} and so we need only to take care of type (a) atoms.
	
Therefore we suppose that $a$ is an $H^p((0,1),dx)$-atom type (a) supported in an interval $I$. If $I\subset \mathcal{J}_j^*$ for some $j\in \mathbb{Z}^*$, then \eqref{eq-proof step 2} follows from part (ii) of Lemma \ref{lem1- atom type ab} . Otherwise, if $I\not\subset \mathcal{J}_j^*$ for any $j\in \mathbb{Z}^*$, then there must exist a largest index $j_1 \in \mathbb{Z}^*$ and a smallest index $j_2 \in \mathbb{Z}^*$  so that $j_1<j_2$ and $I=\sum_{j=j_1}^{j_2}\mathcal{J}_j^*$. 

Set $j_0:=\min\{|j_1|,|j_2|\}$ if $j_1j_2>0$, and $j_0:=0$ if $j_1j_2<0$. Then we have $|I|\sim 2^{-j_0}$. We now decompose $a$ as follows:
		\begin{equation}
	\label{eq1-proof converse}
	a=\sum_{j=j_1}^{j_2} 2^{(j_0-|j|)/p}a_j
	\end{equation}
	where 
	\[
	a_j=2^{-(j_0-|j|)/p}\f{\chi_{\mathcal{J}_j^*}}{\sum_{i\in \mathbb{Z}^*}\chi{\mathcal{J}_i^*}}a.
	\]
	Then $\supp\,a_j\subset \mathcal{J}_j^*$ and
	\begin{equation}
	\label{eq2-proof converse}
	\begin{aligned}
	\|a_j\|_{L^\vc}&\lesi 2^{-(j_0-|j|)/p}\|a\|_{L^\vc}\lesi 2^{-(j_0-|j|)/p}|I|^{-1/p}
	\lesi 2^{-(j_0-|j|)/p}2^{j_0/p}=2^{|j|/p}\sim |\mathcal{J}_j^*|^{-1/p}
	\end{aligned}
	\end{equation}
Therefore, if we write
\[
a_j=\Big[a_j-\f{\chi_{\mathcal{J}_j}}{|\mathcal{J}_j|}\int_{\mathcal{J}_j^*}a_j(x)dx\Big] +\f{\chi_{\mathcal{J}_j}}{|\mathcal{J}_j|}\int_{\mathcal{J}_j^*}a_j(x)dx=a_{j1}+a_{j2}
\]
then from \eqref{eq2-proof converse} it follows that $a_{j1}$ is type (a) atom supported in $\mathcal{J}_j^*$ and that $a_{j2}$ is an type (b) atom. This along with Lemma \ref{lem1- atom type ab} implies that $\|a_j\|_{H^p((0,1),dx)}\lesi 1$. Then taking into account \eqref{eq1-proof converse} we see that $\|a\|_{H^p((0,1),dx)}\lesi 1$, completing our proof.
\end{proof}

\subsection{Fourier--Bessel operators on $((0,1),x^{2\nu+1}dx)$}
Consider the following differential operator
$$ L=-\f{d^2}{dx^2}- \f{2\nu+1}{x}\f{d}{dx}, \qquad \nu>-1.$$

Let $\{\lambda_{k,\nu}\}_{k\ge 1}$ denote the sequence of succesive positive zeros of the Bessel function $J_\nu$ and consider 
$$
\phi_k^\nu(x)=d_{k,\nu}\lambda^{1/2}_{k,\nu}J_\nu(\lambda_{k,\nu}x)x^{-\nu}
$$
where $x\in (0,1)$ and $d_{k,\nu}=\sqrt{2}|\lambda_{k,\nu}J_{\nu+1}(\lambda_{k,\nu})|^{-1}$.

The system $\{\phi_k^\nu\}_k$ forms an orthornomal basis for $L^2((0,1),d\mu)$, where $d\mu(x)=x^{2\nu+1}dx$. It is well known that
\[
L\phi_k^\nu(x)=\lambda_{k,\nu}^2\phi_k^\nu(x).
\]
The operator $L$ has a non-negative self-adjoint extensions which is still denoted by $L$ with domain
\[
D(L)=\{f\in L^2((0,1),x^{2\nu+1}dx): \sum_{k=1}^\infty \lambda_{k,\nu}^4 |\langle f, \phi_k^\nu\rangle|^2<\infty\}.
\]
Let us denote $(X,d,\mu)=((0,1),|\cdot|, x^{2\nu+1}dx)$. One can easily show that
\begin{equation}
\label{eq-measure mu}
\mu(I)=\begin{cases}
x^{2\nu+1}r, \ \ \ &x>r\\
r^{2\nu+2}, \ \ \ &x\le r
\end{cases}
\end{equation}
where $I=(x-r,x+r)\cap (0,\vc)$ with $x\in (0,1)$ and $r< 1$. It is clear then that the triple $(X, |\cdot|, d\mu)$ is a space of homogeneous type with dimension $n=2\nu+2$.

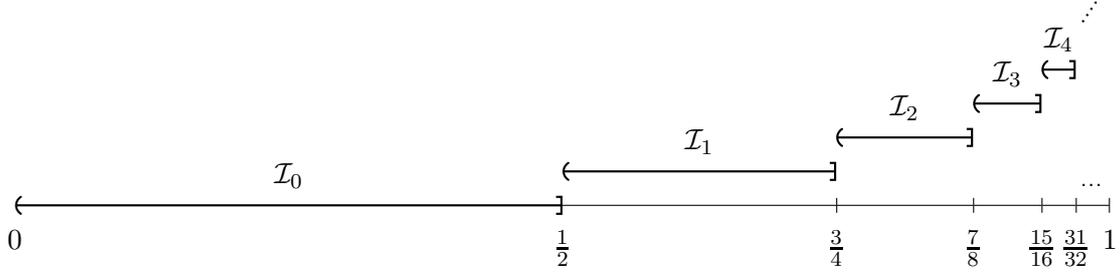
\begin{figure}[h]
\begin{center}
\begin{tikzpicture}[scale=0.9]
\path (-8,-1) coordinate (L5);
\path (0,-1) coordinate (00);
\path (4,-1) coordinate (R1);
\path (6,-1) coordinate (R2);
\path (7,-1) coordinate (R3);
\path (7.5,-1) coordinate (R4);
\path (8,-1) coordinate (R5);

\draw [-] (00)--(R1);
\draw [|-] (R1)--(R2);
\draw [|-] (R2)--(R3);
\draw [|-] (R3)--(R4);
\draw [|-|] (R4)--(R5);

\draw[[-), thick] (00)--(L5);
\draw[[-),thick] (4,-.5) -- (0,-.5);
\draw[[-),thick] (6,0) -- (4,0);
\draw[[-),thick] (7,0.5) -- (6,0.5);
\draw[[-),thick] (7.5,1) -- (7,1);

\node[above=1mm] at (-4,-1) {$\mathcal{I}_{0}$};
\node[above=1mm] at (2,-.5) {$\mathcal{I}_1$};
\node[above=1mm] at (5,0) {$\mathcal{I}_2$};
\node[above=1mm] at (6.5,.5) {$\mathcal{I}_3$};
\node[above=1mm] at (7.25,1) {$\mathcal{I}_4$};
\node[below=2mm] at (00) {$\f{1}{2}$};
\node[below=2mm] at (R1) {$\f{3}{4}$};
\node[below=2mm] at (R2) {$\f{7}{8}$};
\node[below=2mm] at (R3) {$\f{15}{16}$};
\node[below=2mm] at (R4) {$\f{31}{32}$};
\node[below=2mm] at (R5) {$1$};
\node[below=2mm] at (L5) {$0$};

\draw[dotted, thick] (7.6,1.7) -- (7.8,2);
\draw[dotted, thick] (7.6,-.7)--(7.9,-.7);

 \end{tikzpicture}
\end{center}
\caption{Intervals for $\mathcal{I}_j$} \label{fig: dmu}
\end{figure}

As in \cite{DPRS}, we now consider the intervals:
\begin{equation}\label{eq-Ij}
\mathcal{I}_j=(1-2^{-j},1-2^{-j-1}], \ \ \ j=0,1,\ldots.
\end{equation}
which are depicted in Figure \ref{fig: dmu}.
It is obvious that the family $\{\mathcal{I}_j\}_{j\in \mathbb{N}}$ is pairwise disjoint and $(0,1)=\bigcup_{j\in \mathbb{N}}\mathcal{I}_j$. For each $j\in \mathbb{N}$ we shall denote by $\mathcal{I}_j^*=\tfrac{1}{20}\mathcal{I}_j$ and $\mathcal{I}_j^{**}=\tfrac{1}{10}\mathcal{I}_j$.

Consider the following atoms.
\begin{defn}\label{def: atoms 1}
	Let $p\in (\f{2\nu+2}{2\nu+3}, 1]$. A function $a$ is a type (a) atom associated with an interval $I\subseteq (0,1)$ if 
	\begin{enumerate}[\upshape(i)]
		\item $\supp a\subset I$
		\item $\Vert a\Vert_{L^\vc}\le \mu(I)^{-1/p}$
		\item $\int a(x)d\mu(x)=0$
	\end{enumerate}
	A function $a$ is a type (b) atom if
	$$a(x)=\f{\chi_{\mathcal{I}_j}(x)}{\mu(\mathcal{I}_j)^{1/p}}$$ for some $j\in \mathbb{N}$.
	
	We say a function $a$  is  an $H^p((0,1),d\mu)$-atom associated with $I$ if it is either a type (a) or type (b) atom.
\end{defn}

\begin{defn}[Atomic Hardy spaces on $((0,1),\dm)$]
	Given $p\in (\f{2\nu+2}{2\nu+3},1]$, we  say that $f=\sum
	\lambda_ja_j$ is an atomic $H^p((0,1),d\mu)$-representation if
	$\{\lambda_j\}_{j=0}^\infty\in l^p$, each $a_j$ is a $H^p((0,1),d\mu)$-atom,
	and the sum converges in $L^2(X)$. The space $H^p_{at}((0,1),d\mu)$ is then defined as the completion of
	\[
	\left\{f\in L^2:f \ \text{has an atomic
		$H^p((0,1),d\mu)$-representation}\right\},
	\]
	with the norm given by
	$$
	\|f\|_{H^p_{at}((0,1),d\mu)}=\inf\left\{\left(\sum|\lambda_j|^p\right)^{1/p}:
	f=\sum \lambda_ja_j \ \text{is an atomic $H^p((0,1),d\mu)$-representation}\right\}.
	$$
\end{defn}
As in  Definition \eqref{defn-maximal Hardy spaces} we denote by $H^p_{L, \rad}((0,1),d\mu)$ and $H^p_{L, \max}((0,1),d\mu)$ respectively the maximal Hardy spaces defined via radial and non-tangential maximal functions associated to $L$.

Our main result in this section is the following:

\begin{thm}\label{Th1} Let $\nu>-1$.
	
	\noindent  For $p\in (\f{2\nu+2}{2\nu+3},1]$ we have
	\[
	H^p_{L, \rad}((0,1),d\mu)\equiv H^p_{L, \max}((0,1),d\mu)\hookrightarrow H^{p}_{at}((0,1),d\mu).
	\]
	\noindent In the particular case  $p=1$, we have
	\[
	H^1_{L, \rad}((0,1),d\mu)\equiv H^1_{L, \max}((0,1),d\mu)\equiv H^{1}_{at}((0,1),d\mu)
	\]
	with equivalent norms.
\end{thm}
\begin{rem}
It was proved in \cite{DPRS} that  
\[
	H^{1}_{at}((0,1),d\mu)\equiv H^1_{\sqrt{L}, \rad}((0,1),d\mu)\equiv H^1_{\sqrt{L}, \max}((0,1),d\mu)
\]
under the restriction $\nu>-1/2$.  In this situation, Theorem \ref{Th1} gives a maximal function characterization for the Hardy space $H^{1}_{at}((0,1),d\mu)$ for the full range of $\nu$.
\end{rem}

Before we give the proof of Theorem \ref{Th1} we need to develop some technical material. We shall continue to use the notion of ``interval" as expressed in \eqref{eq-interval}. 

We define the \emph{critical function} for $L$ by 
\begin{equation}\label{eq-critical function}
\rho(x) :=\f{1}{3}(1-x), \qquad x\in (1/2,1)
\end{equation}
By a \emph{critical interval} $I_\rho$ we mean an interval with centre $x_I\in (1/2,1)$ and radius $\rho_I:=\rho(x_I)$.


Our critical intervals admit the following desirable properties, whose easy proofs we omit.
\begin{lem}\label{lem-critical functions}
Let $I_\rho$ be a critical interval related to the function $\rho$ from \eqref{eq-critical function}.
	\begin{enumerate}[\upshape(a)]
		\item 	If $y\in I_\rho$ then \;	$ \f{2}{3}\rho_I<\rho(y)<\f{4}{3}\rho_I$.
		\item If $x,y\in I_\rho$ then \;	$ \f{1}{2}\rho(x)<\rho(y)<2\rho(x)$.
	\end{enumerate}
\end{lem}
We emphasize that the restriction $x_I\in (1/2,1)$ for a critical interval is essential in Lemma \ref{lem-critical functions}. It is not difficult to see that these estimates are invalid if we allow $x_I$ to be near the endpoint $0$.

We next have a useful integration by parts lemma.
\begin{lem}\label{IBP}
	Let $f, g$ be sufficiently smooth functions on $X$ and suppose either $f'(1)=0$ or $g(1)=0$ holds. Then we have 
	$$ \int_0^1L f(x)g(x)\,d\mu(x) = \int_0^1 f'(x)g'(x)\,d\mu(x).$$
\end{lem}
\begin{proof}
	From integration by parts and the conditions on $f$ and $g$,
	\begin{align*}
	\int_0^1L f(x)g(x)\,d\mu(x)
	&= -f'(x)g(x)x^{2\nu+1}\Big|_0^1 + \int_0^1 f'(x) (g(x)x^{2\nu+1})' \,dx \\
	&\qquad - \int_0^1 \f{2\nu+1}{x}f'(x) g(x)x^{2\nu+1}\,dx \\
	&= \int_0^1f'(x)g'(x)x^{2\nu+1}\,dx + \int_0^1 f'(x)g(x)(x^{2\nu+1})'\,dx \\
	&\qquad - \int_0^1 \f{2\nu+1}{x}f'(x) g(x)x^{2\nu+1}\,dx \\
	&= \int_0^1f'(x)g'(x)x^{2\nu+1}\,dx.
	\end{align*}
\end{proof}


We now collect together some estimates on the kernels related to $L$. Let $p_t(x,y)$ be the kernel of the heat semigroup $e^{-tL}$ associated with $L$. Then one has the following two sided bounds. 
\begin{lem}[\cite{MSZ} Theorem 1.1]\label{lem-heat kernel}
	For $\nu>-1$ we have 
	$$p_t(x,y) \approx \f{(1+t)^{\nu+2}}{(t+xy)^{\nu+1/2}} \Big(1\land \f{(1-x)(1-y)}{t}\Big)\f{1}{\sqrt{t}}e^{-\f{|x-y|^2}{4t}-\lambda^2_{1,\nu}t} $$
	for $x,y\in X$ and $t>0$.
\end{lem}
Lemma \ref{lem-heat kernel} and a simple calculation yields
\begin{equation}\label{eq-GU pxy}
p_t(x,y)\lesi \f{1}{\mu(I_{\sqrt{t}}(x))}\Big(1\land \f{(1-x)(1-y)}{t}\Big)e^{-\f{|x-y|^2}{ct}}
\end{equation}
for all $x,y\in X$ and $t>0$.

Denote by $q_t(x,y)$ the kernel of $tLe^{-tL}$. It is well-known that the heat kernel bounds \eqref{eq-GU pxy} can be transfered to the kernel $q_t(x,y)$. That is, we have
\begin{equation}\label{eq-GU qxy}
q_t(x,y)\lesi \f{1}{\mu(I_{\sqrt{t}}(x))}e^{-\f{|x-y|^2}{ct}}
\end{equation}
for all $x,y\in X$ and $t>0$.


\begin{lem}\label{lem: heatderiv 1}
	For $\nu>-1$ we have
	\begin{align}\label{heatderiv1}
	|\partial_x p_t(x,y)| \lesssim \f{1}{(xy)^{\nu+1/2}}\f{e^{-(x-y)^2/4t}}{t}
	\end{align} 
	for all $x,y\in X$ and $t>0$.
\end{lem}
\begin{proof}
	Let $\tilde{p}_t(x,y)$ be the heat kernel of the operator Fourier--Bessel operator \eqref{eq-FBdx} on $((0,1),dx)$. It is well known that the heat kernels $\tilde{p}_t(x,y)$ and  $p_t(x,y)$ are related through (see \cite{DPRS})
	$$\tilde{p}_t(x,y)=(xy)^{\nu +1/2}p_t(x,y).$$
It follows that
	\begin{align*}
	&(xy)^{\nu+1/2}\partial_x p_t(x,y) = \big(\partial_x - \tfrac{\nu+1/2}{x}\big) \tilde{p}_t(x,y)
	\end{align*}
and one can now follow the argument in \cite[Lemma 2.4]{DPRS} to obtain
	\[
	\big|\big(\partial_x - \tfrac{\nu+1/2}{x}\big) \tilde{p}_t(x,y) \big| \lesssim \f{e^{-(x-y)^2/4t}}{t}.
	\]
	This yields \eqref{heatderiv1} as desired.
\end{proof}

We now invest Lemma \ref{lem: heatderiv 1} to obtain the following analogue of Lemma \ref{lem: Qest dx}. 
\begin{lem}\label{lem: Qest}
\begin{enumerate}[\upshape (a)]
	\item Let $I$ with $x_I\in (0,1/2)$ and $r_I\ge 1/8$. Then we have
	\begin{align*}
	\Big|\int_{I}q_t(x,y)\dx\Big| \lesssim \f{t^{1/4}e^{-cr_I^2/t}\mu(I)^{1/2}}{y^{\nu+1/2}} \qquad\forall y\in \tfrac{1}{2}I
	\end{align*}
	for any $t>0$. 
	
	\item	For any critical interval $I_\rho$ we have
	\begin{align*}
	\Big|\int_{I_\rho}q_t(x,y)\dx\Big| \lesssim \f{t^{1/4}}{\sqrt{\rho_I}}e^{-c\rho_I^2/t},\qquad \forall y\in\tfrac{1}{2}I_\rho
	\end{align*}
	for any $t>0$. 
	\end{enumerate}
\end{lem}

\begin{proof}
We first prove (a). Choose a cutoff function $\varphi\in C^\infty_c(X)$ supported in $I$ with $0\le \varphi\le 1$, equal to 1 on $\f{3}{4}I$ and with $|\varphi'(x)|\lesssim 1/r_I$.
	Then 
	\begin{align*}
	\Big|\int_{I}q_t(x,y)\dx\Big|&\lesi  \Big|\int_{X}t\partial_tp_t(x,y)\varphi(x)\dx\Big|+\Big|\int_{I\backslash \tfrac{3}{4}I }q_t(x,y)[1-\varphi(x)]\dx\Big|
	=:I_1+I_2.
	\end{align*}
	For the term $I_2$, using \eqref{eq-GU qxy} and that $|x-y|\sim r_I$, we have
	\[
	\begin{aligned}
	I_2
	\lesi \int_{I\backslash \tfrac{3}{4}I }\f{1}{\mu(I_{\sqrt{t}}(x))}e^{-\f{|x-y|^2}{ct}}\dx
	\lesi e^{-\f{r_I^2}{c't}}\int_{I\backslash \tfrac{3}{4}I }\f{1}{\mu(I_{\sqrt{t}}(x))}e^{-\f{|x-y|^2}{2ct}}\dx 	
	\lesi e^{-\f{r_I^2}{c't}}.
	\end{aligned}
	\]
For the term $I_1$ since $\partial_tp_t(\cdot,y)=-L p_t(\cdot,y)$, then by Lemma \ref{IBP}, estimate \eqref{heatderiv1} and the fact that $|\varphi'(x)|\lesi r_I^{-1}$, we obtain
	\[
	\begin{aligned}
	I_1
	=\Big|\int_{X}t\partial_{x}p_t(x,y)\partial_{x}\varphi(x)\dx\Big|
	\lesi \f{1}{r_I}\Big|\int_{\f{3}{4}I}\f{1}{(xy)^{\nu+1/2}}e^{-(x-y)^2/ct}\dx\Big|.
	\end{aligned}
	\]
	H\"older's inequality and the estimate $|x-y|\sim r_I$ then gives
	\begin{equation}
	\label{eq-int qtxy bessel 2}
		\begin{aligned}
	I_{1}
	&\lesi \f{1}{r_I} \mu(I)^{1/2}\Big[\int_{\f{3}{4}I}\f{1}{(xy)^{2\nu+1}}e^{-(x-y)^2/c't}\,\dm(x)\Big]^{1/2}\\
	&\lesi \f{t^{1/4}}{r_I} \mu(I)^{1/2}e^{-cr_I^2/t}\Big[\int_{\f{3}{4}I}\f{1}{y^{2\nu+1}}\f{1}{\sqrt{t}}e^{-(x-y)^2/2c't}dx\Big]^{1/2} \\
	&\lesi \f{t^{1/4}}{r_I} \f{\mu(I)^{1/2}e^{-cr_I^2/t}}{y^{\nu+1/2}}.
	\end{aligned}
	\end{equation}
	The condition $r_I\ge 1/8$ allows us to conclude
	\[
	I_{1}\lesi \f{t^{1/4}e^{-cr_I^2/t}\mu(I)^{1/2}}{y^{\nu+1/2}}.
	\]
	This and the estimate of $I_2$ proves (a).
	
	If $I_\rho$ is a critical interval, then we have $y\sim x_I$ and $\mu(I)\sim x_I^{2\nu+1}\rho_I$. Hence, \eqref{eq-int qtxy bessel 2} implies
	\[
	I_1\lesi \f{t^{1/4}}{\sqrt{\rho_I}}e^{-c\rho_I^2/t}.
	\]
	Part (b) now follows by combining this estimate  with the estimate of $I_2$ from (a).	
\end{proof}

We can now give the analogue of Lemma \ref{lem1- atom type ab} for $p=1$. Note that the intervals $\II^*_j$ and $\II^{**}_j$ have been defined in the comments after \eqref{eq-Ij}.
\begin{lem}\label{lem2- atom type ab}
Let that $\nu>-1$.  Suppose that $a$ is either 
\begin{enumerate}[\upshape(i)]
	\item  a type-(b) $H^1\big((0,1),\dm\big)$-atom, or 
	\item  a type-(a) $H^1((0,1),\dm)$-atom supported in $\mathcal{I}_j^*$ for some $j\in \mathbb{N}$.
	\end{enumerate}
Then there exists $C>0$ independent of $a$ so that
	\[
	\Big\|\sup_{t>0}|e^{-tL}a|\Big\|_{L^1((0,1),\dm)}\le C.
	\]
\end{lem}

\begin{proof}[Proof of Lemma \ref{lem2- atom type ab}]
Part (i). 
	Suppose that $a$ is an $H^1((0,1),d\mu)$-atom of type-(b). Then $a=\mu(\mathcal{I}_j)^{-1}\chi_{\mathcal{I}_j}$ for $j\ge 0$. Since the radial maximal operator associated with $L$ is uniformly bounded on $L^\infty$ then we have, for $j=0,1,2$,
	\[
	\Big\|\sup_{t>0}|e^{-tL}a|\Big\|_{L^1((0,1),d\mu)}\le \|a\|_{L^\vc}=\mu(\mathcal{I}_j)\lesi 1.
	\]
	In a similar way we have, for $j\ge 3$, 
	\[
	\Big\|\sup_{t>0}|e^{-tL}a|\Big\|_{L^1(3\mathcal{I}_j,\dm)}\lesi 1,
	\]
	and so it remains to show
	\begin{align}\label{eq-atom type ab1}
	\Big\|\sup_{t>0}|e^{-tL}a|\Big\|_{L^1((3\mathcal{I}_j)^c,\dm)}\lesi 1.
	\end{align}
	Firstly note that if $y\in \mathcal{I}_j$ then we have the estimates $y\sim 1$ and $1-y\sim r_{\mathcal{I}_j}\sim \mu(\mathcal{I}_j)\sim 2^{-j}$. Secondly if $x\in (3\mathcal{I}_j)^c$ then we also have $|x-y|\sim 1-x$. These facts along with Lemma \ref{lem-heat kernel} gives
    $$
    \begin{aligned}
    \Big\|\sup_{t>0}|e^{-tL}a|\Big\|_{L^1((3\mathcal{I}_j)^c,\dm)}&\lesi \mu(\mathcal{I}_j)^{-1}\int_{(3\mathcal{I}_j)^c}\sup_{t>0}\int_{\mathcal{I}_j} \f{(1-x)}{(t+x)^{\nu+1/2}}\f{2^{-j}}{t^{3/2}}e^{-\f{(1-x)^2}{ct}} dy\, \dx\\
    &\lesi \int_{(3\mathcal{I}_j)^c} \f{2^{-j}x^{\nu+1/2}}{(1-x)^2} dx\\
    \end{aligned}
	$$
	Now since $(3\mathcal{I}_j)^c=(0,1-3\x 2^{-j})$ then a direct calculation gives \eqref{eq-atom type ab1}.



Part (ii).
	Suppose that $a$ is an $H^1((0,1),d\mu)$-atom of type-(a) associated to an interval $I\subset \mathcal{I}_j^*$ for some $j\in \mathbb{N}$. 
		
		We now consider two cases: $j\ge 1$ and $j=0$.
		
	\textbf{Case 1: $j\ge 1$.} 
	In a similar fashion to part (ii)  of Lemma \ref{lem1- atom type ab}, it suffices to prove that 
	\begin{align}\label{eq-atom type ab2}
	\Big\|\sup_{t>0}|e^{-tL}a|\Big\|_{L^1((2I)^c,\dm)}\lesi 1.
	\end{align}
	
	In this situation we have $x\sim y \sim 1$ for every $x\in 3\mathcal{I}_j$ and $y\in I$. This fact, the cancellation property of $a$, and Lemma \ref{lem: heatderiv 1} imply that
	\begin{equation}\label{eq1- lem2 atom a}
	\begin{aligned}
	 &\Big\|\sup_{t>0}|e^{-tL}a|\Big\|_{L^1((3\mathcal{I}_j\backslash 2I),\dm)}\\
	 &\qquad=\int_{3\mathcal{I}_j\backslash 2I}\sup_{t>0}\Big|\int_I [p_t(x,y)-p_t(x,x_I)]a(y)\dy\Big|\dx\\
	 &\qquad \lesi \int_{3\mathcal{I}_j\backslash 2I}\sup_{t>0}\int_I  \f{|y-x_I|}{(xy)^{\nu+1/2}}\f{e^{-(x-y)^2/ct}}{t}|a(y)|\,\dy\,\dx.
	\end{aligned}
	\end{equation}
	Since $d\mu(x)\sim dx$ whenever $x\in 3\mathcal{I}_j$, and $|x-y|\sim |x-x_I|$ whenever $x\in (2I)^c$ and $y\in I$, we may continue with
	\begin{equation*}
	\begin{aligned}
\Big\|\sup_{t>0}|e^{-tL}a|\Big\|_{L^1((3\mathcal{I}_j\backslash 2I),\dm)}
	  \lesi \int_{3\mathcal{I}_j\backslash 2I}\sup_{t>0}\int_I  \f{r_I}{t}e^{-|x-x_I|^2/ct}|a(y)|\dy dx
	\lesi \int_{3\mathcal{I}_j\backslash 2I}\f{r_I}{|x-x_I|^2}dx\\
	\end{aligned}
	\end{equation*}
which evaluates to
    \begin{equation}\label{eq1-proof lem 2 atom (a)}
    \Big\|\sup_{t>0}|e^{-tL}a|\Big\|_{L^1(3\mathcal{I}_j\backslash 2I)}\lesi 1.	
    \end{equation}

Since $y\sim 1$ for $y\in I\subset 3\mathcal{I}_j$ then arguing similarly to \eqref{eq1- lem2 atom a} we have
    \begin{equation*}
    \begin{aligned}
\Big\|\sup_{t>0}|e^{-tL}a|\Big\|_{L^1((3\mathcal{I}_j)^c,\dm)}
    & \lesi \int_{(3\mathcal{I}_j)^c}\sup_{t>0}\int_I  \f{|y-x_I|}{(xy)^{\nu+1/2}}\f{1}{t}e^{-(x-y)^2/ct}|a(y)|\,\dy\,\dx\\
    & \lesi \int_{(3\mathcal{I}_j)^c}\sup_{t>0}\int_I  \f{r_I}{x^{\nu+1/2}}\f{1}{t}e^{-|x-x_I|^2/ct}|a(y)|\,\dy\,\dx\\
    & \lesi \int_{(3\mathcal{I}_j)^c}\f{r_I}{|x-x_I|^2} x^{\nu+1/2}dx.
    \end{aligned}
    \end{equation*}
Noting that $(3\mathcal{I}_j)^c=(0,1-3\x 2^{-j})$, and applying the estimates $r_I\lesi 2^{-j}$ and $|x-x_I|\sim |1-x|$ whenever $x\in (3\mathcal{I}_j)^c$ we continue with
    \begin{align}\label{eq-atom type ab3}
\Big\|\sup_{t>0}|e^{-tL}a|\Big\|_{L^1((3\mathcal{I}_j)^c,\dm)}
     & \lesi \int_{0}^{1-3\x 2^{-j}}\f{2^{-j}x^{\nu+1/2}}{|1-x|^2} dx\lesi 1.
    \end{align}
    Estimate \eqref{eq-atom type ab3} and \eqref{eq1-proof lem 2 atom (a)} together gives \eqref{eq-atom type ab2}.
 
 \textbf{Case 2: $j= 0$.} We write
 \[
 \begin{aligned}
 &\Big\|\sup_{t>0}|e^{-tL}a|\Big\|_{L^1((0,1),d\mu)}\\
 &\qquad \le \Big\|\sup_{t\ge 1}|e^{-tL}a|\Big\|_{L^1((0,1),d\mu)}+\Big\|\sup_{0<t<1}|e^{-tL}a|\Big\|_{L^1(\mathcal{I}_0^{**})}
  +\Big\|\sup_{0<t<1}|e^{-tL}a|\Big\|_{L^1((\mathcal{I}_0^{**})^c)}\\
  &\qquad=:A_1+A_2+A_3.
 \end{aligned}
 \]
 From Lemma \ref{lem-heat kernel} we see that the heat kernel is dominated by 1 whenever $t\ge 1$, so as a consequence we obtain easily that $A_1\lesi \|a\|_{L^1}\lesi 1$.
 Again from Lemma \ref{lem-heat kernel} and the fact that $|x-y|\sim x \sim 1$ we have
 \[
 \begin{aligned}
 A_3\lesi \int_{(\mathcal{I}_0^{**})^c}\sup_{0<t<1}\int_I \f{1}{(t+y)^{\nu+1/2}} \f{1}{\sqrt{t}}e^{-\f{c}{t}}|a(y)|\dy \dx. 
 \end{aligned}
 \]
 Since the inner integrand can be controlled by a constant multiple of $|a(y)|$ 
 for all $t\in (0,1)$ and $y\in I$, then we obtain $A_3\lesi \|a\|_{L^1}\lesi 1$.
 
To handle  $A_2$  we shall employ a comparison with a Bessel operator  on $((0,\vc), d\mu)$ defined by
 \[
 \mathfrak{L}=-\f{d^2}{dx^2}-\f{2\nu+1}{x}\f{d}{dx}.
 \]
Let $h_t(x,y)$ be the kernel of $e^{-t\mathfrak{L}}$. Then it is well-known that
\begin{align}
\label{eq-fullbessel hk1}
 h_t(x,y) &\lesi \f{1}{\mu(I_{\sqrt{t}}(x))}e^{-\f{|x-y|^2}{ct}}  \\
\label{eq-fullbessel hk2}
  |\partial_x h_t(x,y)| &\lesi \f{1}{\sqrt{t}\,\mu(I_{\sqrt{t}}(x))}e^{-\f{|x-y|^2}{ct}}
\end{align}
 We now split the term $A_2$ as follows:
 \[
 \begin{aligned}
 A_2&\le \Big\|\sup_{0<t<1}|(e^{-tL}-e^{-t\mathfrak{L}})a|\Big\|_{L^1(\mathcal{I}_0^{**})}+\Big\|\sup_{0<t<1}|e^{-t\mathfrak{L}}a|\Big\|_{L^1(\mathcal{I}_0^{**})}
 =:A_{21}+A_{22}.
 \end{aligned}
 \]
 From \eqref{eq-fullbessel hk1} and \eqref{eq-fullbessel hk2} above, by a standard argument we have 
 \begin{align*}
 A_{22}\le \Big\|\sup_{t>0}|e^{-t\mathfrak{L}}a|\Big\|_{L^1(0,\vc)} \lesssim 1.
 \end{align*}
 On the other hand,  Corollary 4.7 in \cite{DPRS} shows that $A_{21}\lesi \|a\|_{L^1}\lesi 1$ and hence we have $A_2\lesi 1$.
 
 Taking the estimates of $A_1, A_2, A_3$ into account, we arrive at the required estimate as stated in the Lemma.
This completes our proof.
\end{proof}


We are now ready to prove Theorem \ref{Th1}, which is a direct consequence of Theorem \ref{mainthm1} and the following proposition.
\begin{prop}
	\label{prop-Thm1 bessel 2}
	Let $\nu>-1$. Then for $p\in (\f{2\nu+2}{2\nu+3},1]$ we have
	\begin{enumerate}[{\rm (i)}]
		\item $H^p_{L,\rad}(X)\subset H^p_{at}(X)$;
		\item $H^1_{at}(X)\subset H^1_{L,\rad}(X)$.
	\end{enumerate}
\end{prop}

\begin{proof} Part (i). Suppose that $a$ is a $(p,M)_L$-atom as in Definition \ref{def: L-atom} associated to an interval $I$. We consider two cases: $4I\cap [1,\vc)\ne \emptyset$ and $4I\subset (0,1)$. 
	
\noindent\textbf{Case 1: $4I\cap [1,\vc)\ne \emptyset$.} In this situation, it easy to see that if $x_I\in \mathcal{I}_j$ for some $j\ge 0$, then $I\subset \mathcal{I}_{j-1}\cup \mathcal{I}_j\cup \mathcal{I}_{j+1}$ and $\mu(I)\sim \mu(\mathcal{I}_j)$. Hence, using the decomposition
\[
a=\Big[a-\f{\chi_{\mathcal{I}_j}}{\mu(\mathcal{I}_j)}\int a\Big] + \f{\chi_{\mathcal{I}_j}}{\mu(\mathcal{I}_j)}\int a=:\tilde{a}_1+\tilde{a}_2
\]
we see that $a_1$ is a type-(a) atom, while $a_2$ is a type-(b) atom from Definition \ref{def: atoms 1}. Hence, $a\in H^p_{at}((0,1),d\mu)$.

\medskip

\noindent\textbf{Case 2: $4I\subset (0,1)$.} In this case, $a$ can be expressed in the form $a=Lb$. We now write
	\[
	a=Le^{-r_I^2L}b+L(I-e^{-r_I^2L})b=Le^{-r_I^2L}b+(I-e^{-r_I^2L})a=a_1+a_2
	\]
where $b$ is supported in $B$ and satisfies
\[
\|b\|_{L^\vc}\le r_I^2 \mu(I)^{-1/p}.
\]
At this stage the proof is similar to that of Theorem \ref{mainth-Dirichlet} and we will just sketch the main ideas.

We take care  of $a_2$ only, since $a_1$ can be treated similarly. We  choose $k_0\in \mathbb{N}$ so that $2^{k_0-1}r_I\le 1-x_I <2^{k_0}r_I$. Hence $k_0\ge 3$. We set $S_j(I)=[2^{j+1}I\backslash 2^jI]\cap (0,1)$ if $j>0$ and $S_0(I)=2I$. Then as in  the proof of Theorem \ref{mainth-Dirichlet} we decompose $a_2$ as follows:
	\[
	\begin{aligned}
	a_2=\ &\sum_{j=k_0-3}^\vc a_2\chi_{S_j(I)}+\sum_{j=0}^{k_0-3}\Big(a_2\chi_{S_j(I)}-\f{\chi_{S_j(I)}}{\mu(S_j(I))}\int_{S_j(I)}a_2\Big)\\
	&+\sum_{j=0}^{k_0-3}\Big(\f{\chi_{S_j(I)}}{\mu(S_j(I))}-\f{\chi_{S_{j-1}(I)}}{\mu(S_{j-1}(I))}\Big)\int_{2^{k_0-3}\backslash 2^jI}a_2+\f{\chi_{2I}}{\mu(2I)}\int_{2^{k_0-3}I}a_2\\
	=:\ & A_1 + A_2 + A_3 + A_4.
	\end{aligned}
	\]
	Arguing similarly to the proof of Theorem \ref{mainth-Dirichlet} we can show that $A_1$, $A_2$ and $A_3$ can each be expressed as atomic representations of atoms from Definition \ref{def: atoms 1}; $A_1$ as type-(b) atoms,  $A_2$ and $A_3$ as type-(a) atoms.
%

	It remains to take care of $A_4$. We consider two subcases: $x_I\le 1/2$ and $x_I>1/2$.
	
	\textbf{Subcase 2.1: $x_I\in (0,1/2]$.} In this situation we have $2^{k_0}I=(0,1)$ and hence $2^{k_0-3}r_I\ge 1/8$. We now have
	\[
	\begin{aligned}
	\int_{2^{k_0-3}I}a_2
	&=\int_0^{r_I^2}\int_I \int_{2^{k_0-3}I}q_s(x,y)a(y)\dx\dy \f{ds}{s}\\
	&\lesi \mu(2^{k_0-3}I)^{1/2}\int_0^{r_I^2}s^{1/4}e^{-c2^{2k_0}r_I^2/s}\int_I \f{a(y)}{y^{\nu+1/2}}\dy \f{ds}{s}\\
	 &\lesi e^{-c2^{2k_0}}\mu(2^{k_0-3}I)^{1/2}\int_I \f{a(y)}{y^{\nu+1/2}}\dy
	\end{aligned}
	\]
	where in the second inequality we used (a) in Lemma \ref{lem: Qest}.
	By  H\"older's inequality we have
	\begin{align*}
	\int_I \f{a(y)}{y^{\nu+1/2}}\dy
	\le  \mu(I)^{1/2} \Big[\int_I \f{|a(y)|^2}{y^{2\nu+1}}\dy\Big]^{1/2}
	\lesssim  \mu(I)^{1/2-1/p}
	\end{align*}
and so we may continue with 
	 \[
	 \begin{aligned}
	 \int_{2^{k_0-3}I}a_2& \lesi e^{-c2^{2k_0}}\mu(2^{k_0-3}I)^{1/2} \mu(I)^{1/2-1/p}.
	 \end{aligned}
	 \]
	Now since $2^{k_0}I=(0,1)$, then $\mu(2^{k_0-3}I)\sim 1 \sim \mu(\mathcal{I}_0)$, and so we have
	\[
	\|A_4\|_{L^\vc}\lesi \f{e^{-c2^{2k_0}}\mu(2^{k_0-3}I)^{1/2} \mu(I)^{1/2-1/p}}{\mu(I)} \lesi \f{1}{\mu(2^{k_0-3}I)^{1/p}}\sim \f{1}{\mu(\mathcal{I}_0)^{1/p}}.
	\]
%
	We now decompose $A_4$ as follows:
	\[
	A_4= \Big[A_4-\f{\chi_{\mathcal{I}_0}}{\mu(\mathcal{I}_0)}\int A_4\Big] + \f{\chi_{\mathcal{I}_0}}{\mu(\mathcal{I}_0)}\int A_4=:A_{41}+A_{42}.
	\]
	Then $A_{41}$ is an $H^p((0,1),d\mu)$ of type-(a) and $A_{42}$ is an $H^p_{at}((0,1),d\mu)$-atom of type-(b). We conclude therefore that $\|A_{4}\|_{H^p_{at}((0,1),d\mu)}\lesi 1$.
	
	\textbf{Subcase 2.2: $x\in (1/2,1)$.} We may repeat the argument in Subcase 2.1 but make use of Lemma \ref{lem: Qest} (b) rather than Lemma \ref{lem: Qest} (a). 
	This completes the proof of part (i).
	
	\medskip

	Part (ii).  The argument from Step 2 of Theorem \ref{Th1-Bessel 1} carries over harmlessly but with Lemma \ref{lem2- atom type ab} in place of Lemma \ref{lem1- atom type ab}. We leave the details to the interested reader.
\end{proof}

{\bf Acknowledgement.} T. A. Bui and X. T Duong were supported by the research grant ARC DP140100649 from the Australian Research Council.


\end{document}